\documentclass{article}
\usepackage{xcolor}
\usepackage{amsfonts}
\usepackage{amsmath, amsfonts, amsthm, amssymb, amscd}
\usepackage[titletoc]{appendix}
\usepackage[mathcal]{euscript}
\usepackage{mathrsfs}
\usepackage{dsfont}
\usepackage{ulem}
\usepackage{bbm}
\usepackage{verbatim}
\usepackage{slashed}
\usepackage{etoolbox}
\usepackage[hidelinks]{hyperref}
\usepackage{titlesec}
\usepackage{a4wide}
\usepackage{hyperref} 

\usepackage{lineno} 

\newtheorem{theorem}{Theorem}[section]
\newtheorem{proposition}[theorem]{Proposition}
\newtheorem{lemma}[theorem]{Lemma}
\newtheorem{corollary}[theorem]{Corollary}

\newtheorem{remark}[theorem]{Remark}

\numberwithin{equation}{section}

\usepackage{multirow}
\usepackage{tabularx}
\usepackage{lscape}

\newcommand \Kcal {\mathcal K}

\newcommand \Hcal {\mathcal H}

\newcommand \Lcal {\mathcal L}

\newcommand \delb {\bar {\del}}

\newcommand \Hb{\bar H}

\newcommand \etau{\underline{\eta}}

\newcommand \del \partial
\newcommand \delu {\uline{\del}}

\newcommand \Tu {\uline{T}}

\newcommand \hu {\uline{h}}

\newcommand \Au{\uline{A}}

\newcommand \Psiu{\uline{\Psi}}
\newcommand \Phiu{\uline{\Phi}}

\newcommand \Abf {{\bf A}}
\newcommand \Bbf {{\bf B}}

\newcommand{\rank}{\text{rank}}
\newcommand{\ord}{\text{ord}}
\newcommand{\con}{\mathrm{con}}

\newcommand \RR{\mathbb{R}}
\newcommand \eps{\varepsilon}

\newcommand {\vep}{\varepsilon}

\newcommand {\Ebf}{{\bf E}}
\newcommand {\ebf}{{\bf e}}

\AtBeginEnvironment{appendices}{%
	\addtocontents{toc}{\protect} 
}

\title{A non-linear damping structure and global stability of wave-Klein-Gordon coupled system in $\RR^{3+1}$}

\author{
Yue Ma\textsuperscript {1,2},\hskip.2cm
Weidong Zhang\textsuperscript{3,4,5}
}

\begin{document}

\maketitle
\footnotetext[1]{School of Mathematics and Statistics, Xi'an Jiaotong University, Xi'an, 710049 Shaanxi, People’s Republic of China. E-mail: {\tt yuemath@xjtu.edu.cn}}
\footnotetext[2]{Partially supported by the Fundamental Research Funds of Xi'an Jiaotong University (Grant No.:xzy012023034) and Shaanxi Fundamental Science Research Project for Mathematics and Physics (Grant No. 23JSY032)}
\footnotetext[3]{School of Mathematics and Statistics, Xi'an Jiaotong University, Xi'an, 710049 Shaanxi, People’s Republic of China. E-mail: {\tt zwd13892650621@stu.xjtu.edu.cn}}
\footnotetext[4]{Laboratoire Jacques-Louis Lions, Sorbonne Universit\'e, 75252 Paris, France. Email: {\tt weidong.zhang@sorbonne-universite.fr}}
\footnotetext[5]{Partially supported by the CSC scholarship program (Project ID: 202406280362) and the Fundamental Research Funds of Xi'an Jiaotong University (Grant No.:xzy022023006).}

\begin{abstract}
	This paper establishes the global existence of solutions for a class of wave-Klein-Gordon coupled systems with specific nonlinearities in 3+1-dimensional Minkowski spacetime. The study demonstrates that imposing certain constraints on the coefficients of these specific nonlinear terms induces a damping effect within the system, which is crucial for proving the global existence of solutions. The proof is conducted within the framework of a bootstrap argument, primarily employing the hyperboloidal foliation method and the vector field method. 
\end{abstract}

	\section{Introduction}
	
	We investigate the global existence for a class of coupled systems of nonlinear wave-Klein-Gordon equations in (3+1) Minkowski space-time dimensions:
	\begin{equation}\label{new1.1}
		\aligned
		&\Box u = A^{\alpha\beta}\del_{\alpha}v\del_{\beta}v + Bv^2,
		\\
		&\Box v + c^2v = H^{\alpha\beta}u\del_\alpha\del_\beta v + P^\alpha u\del_\alpha v + Quv,
		\endaligned
	\end{equation}
	with initial data:
	\begin{equation}\label{eq02}
		u\vert_{t=2} = u_0,\quad \del_t u\vert_{t=2} = u_1,\quad v\vert_{t=2} = v_0,\quad \del_t v\vert_{t=2} = v_1
	\end{equation}
	supported in the disc $\{|x|<1\}$, and satisfying the smallness condition
	\begin{equation}\label{eq1-24-0march-2025}
	\|u_0\|_{H^{N+1}(\RR^3)} + \|u_1\|_{H^N(\RR^3)} + \|v_0\|_{H^{N+1}(\RR^3)} + \|v_1\|_{H^N(\RR^3)}\leq \eps.
	\end{equation}
	Here $A^{\alpha\beta}$, $B$ and $c$ are constants, $c>0$. The Minkowski metric $\eta = dt^2 - \sum_{a=1}^3(dx^a)^2$.

When $P^\alpha=Q=0$, this type of system was introduced in \cite{Lefloch-Ma-2016-CMP} (see also \cite{Wang-2020}) as a simplified model of the Einstein-Klein-Gordon system, which captures the main difficulties of the full physical system. It was then studied extensively in mathematical literature, such as \cite{Ionescu-Pausader-2014,Chen-Lindblad-2023,Duan-Ma-Zhang-2023}. Later, S. Dong and Z. Wyatt initiate the study for the case $Q\neq 0$ in \cite{Dong-Wyatt-2020}, see also \cite{Chen-Lindblad-2023,Xinyu Cheng-2023,Ionescu-Pausader-2019,Zhimeng Ouyang-2023,Qian Zhang-2024,Ozawa-Tsutaya-Tsutsumi-1995,Tsutsumi-2003}.

In this article, we will make a first glance on the case $P\neq 0$, which is essentially different from the two previous ones. To be more precises, let us consider the following three ODEs:
\begin{subequations}
	\begin{align}
		&(1+H)v'' + v =0,\label{eq1-29-08-2024}
		\\
		&v'' + Pv' + v = 0,\label{eq2-29-08-2024}
		\\
		&v'' +(1+Q)v = 0.\label{eq3-29-08-2024}
	\end{align}
\end{subequations}
It is clear that the first and last one, providing that $H$ and $Q$ being sufficiently small, always have their characteristic roots purely imaginary which leads to an effect of oscillation, while the second one always has characteristic roots with non-trivial real part, which leads to exponentially growth or damping, depending on its sign. Similar story occurs in the PDE case of \eqref{new1.1}, where the sign of the Minkowski inner product $(P^{\alpha}u\del_{\alpha},S)_{\eta}$ takes the role of the sign of $P$ in \eqref{eq2-29-08-2024}. Here $S = t\del_t +r\del_r$. Roughly speaking, we expect $(P u,S)_{\eta}\geqslant 0$ in order to avoid the experiential growth.

In the present work, in order to concentrate on the essential point, we regard the case where $H = Q = 0$ and $P = (1,0,0,0)$:
\begin{equation}\label{eq01}
	\aligned
	&\Box u = A^{\alpha\beta}\del_{\alpha}v\del_{\beta}v + Bv^2,
	\\
	&\Box v + c^2v = u\del_t v.
	\endaligned
\end{equation}
The technical developed below can be easily applied to the case with general $P$. Then combining with \cite{Lefloch-Ma-2016-CMP,Chen-Lindblad-2023} etc., one can analysis the full system \eqref{new1.1}.

To specify the damping condition, we need to introduce

\begin{equation}\label{eq5-29-08-2024}
\underline{A}^{00} = A^{00} - 2\frac{x^a}{t}A^{a0} + \frac{x^ax^b}{t^2}A^{ab},
\end{equation}
and the nonlinear damping condition is written as
	\begin{equation}\label{eq1-24-07-2024}
		\underline{A}^{00} +
		\frac{B}{c^2}(s/t)^2\leqslant 0.
	\end{equation}	
Two simple examples satisfying the above condition are: 
	\begin{itemize}
		\item $A = 0$, $B<0$,
		\item $A$ is negative defined and $\vert B\vert/c^2$ is controlled by the smallest eigenvalue of $A$.
	\end{itemize}
	 We emphasize that the quadratic form $A$ need not to be null.
	 
This condition \eqref{eq1-24-07-2024} seems to be artificial, however, it is closely linked to the {\sl energies conditions} in general relativity. As we will see, the essence of  \eqref{eq1-24-07-2024} is a negative quadratic contribution from the Klein-Gordon component on wave source (modulo high-order and faster decaying terms). If we turn back to the complete Einstein-Klein-Gordon system (from where \eqref{eq01} is modeled), this negative contribution forces the light-cones bending in interior, which coincides with the physical intuition that the mass (Klein-Gordon field) forces the light rays bend towards the mass center.


In the present article, for the clarity of the main technique, we restrict our discussion to the compactly support case, i.e. we suppose that the initial data $u_\ell$, $v_\ell$, $\ell = 0, 1$ are compactly supported in $\mathcal{H}_2^* = \mathcal{H}_2\cap\mathcal{K}$ with $\mathcal{H}_2 = \{(t,x)\vert t = \sqrt{s^2 + \vert x\vert^2}\}$, $\mathcal{K} = \{r < t - 1\}$ and $\mathcal{H}_2^* = \mathcal{H}_2\cap \mathcal{K}$. This is not an essential restriction. In fact one can apply the Euclidean-hyperboloidal foliation, and treat the non-compactly supported initial data with sufficient decay rate at spatial infinity.

	
	
The global well-posedness of nonlinear wave and Klein-Gordon equations in $3+1$ space-time was initialized by by S.Klainerman, D. Christodoulou and J. Shatah, see \cite{Klainerman-1980,Klainerman-1985}, \cite{Christodoulou-1986} and \cite{Shatah-1985}. It has come to our attention that the question of the coupled wave-Klein-Gordon-type systems' (or systems with similar nature) global existence is equally a matter of significant interest. A pioneer work \cite{Bachelot-1988} of A. Bachelot solved the global existence problem for the Dirac-Klein-Gordon system in $\RR^{3+1}$-dimensional Minkowski space-time with massless scalar field coupled with massive Dirac one, and proved the sharp pointwise decay and scattering result. Psarelli \cite{Psarelli-2005} proved the global existence and asymptotic estimate for the small amplitude solutions to the massive Maxwell–Dirac system in the 3+1-dimensional Minkowski space-time. It is particularly worth mentioning that, S. Klainerman, Q. Wang, and S. Yang presented a method of global existence for the Maxwell-Klein-Gordon system with non-compact support initial data in \cite{Klainerman-Wang-Yang-2020}. A. Ionescu and B. Pausader \cite{Ionescu-Pausader-2022}, P. LeFloch and Y. Ma \cite{LeFloch-Ma-2024} established the global existence of the Einstein-Klein-Gordon system with non-compactly supported initial data. In \cite{Stingo-2015}, A. Stingo proved the global stability of a class of quasi-linear Klein-Gordon equation in one-space dimension with non-compactly supported initial data. In \cite{Ifrim-Stingo-2019}, A. Stingo and M. Ifrim proved almost global well-posedness for quasilinear strongly coupled wave-Klein-Gordon systems with small and localized data in two space dimensions. For other recent results on wave-Klein-Gordon systems, see \cite{Lindblad-Rodnianski-2003,Lindblad-Rodnianski-2010,MinggangCHENG-2021,Shijie Dong-2021,Shijie Dong-2024,Dong-Ma-Yuan,Dong-Wyatt-2024,Georgiev-1990,Georgiev-1991,Katayama-2012,Katayama-Matsumura-Sunagawa-2015,Stingo-2023,Qian Zhang-2025,Yang-2015,Andersson-Bar-2018,Boussaid-Comech-2018,Yang-Yu-2019,Selberg-Tesfahun-2021,Nolasco-2021,Dai-Mei-Wei-Yang-2025,Cacciafesta-Danesi-Meng-2024,Matsuno-Ueno-2023,Herr-Ifrim-Spitz-2024,Miao-Yang-Yu-2025}.

In the present article we focus on the mixed quadratic terms coupled in Klein-Gordon equation, which contain both wave and Klein-Gordon components. Most of these terms are trivial in 3+1 dimensional space-time \cite{Lefloch-Ma-2014}, while three are more delicate : $u\del_{tt} v$, $u\del_t v$, and $uv$. Among them, $u\del_{tt} v$ was solved in \cite{Lefloch-Ma-2016-CMP}, while $uv$ was handled by the work of S. Dong and Z. Wyatt \cite{Dong-Wyatt-2020}. The present result will be the final piece of the puzzle to complete the second-order wave-Klein-Gordon interaction in 3+1 dimensional space-time. 
	 
Now we state the main result.
	
\begin{theorem}\label{theorem 1.1}
Consider the wave-Klein-Gordon coupled system \eqref{eq01} with initial data \eqref{eq02} and the damping condition \eqref{eq1-24-07-2024}. Then for any positive integer $N \geqslant 7$, there exists $\varepsilon_0 > 0$ such that if \eqref{eq1-24-0march-2025} is verified with $0\leqslant \varepsilon \leqslant \varepsilon_0$, then the corresponding local solution extends to time infinity.
\end{theorem}
	
	The structure of this article is as follows. 
	In Sections 2 and 3, we recall the technical ingredients in the hyperboloidal framework. Sections 4 - 6 are devoted to the proof of Theorem \ref{theorem 1.1}, which is composed by three steps. In the first step we rely on the Klainerman-Sobolev inequalities \eqref{eq3-21-04-2023}, \eqref{eq5-14-08-2023} and the bootstrap assumption \eqref{Energy bound}, and obtain preliminary $L^2$ and $L^\infty$ estimations. In the second step we rely on the linear estimates Propositions~\ref{prop-01-10-2023}, Proposition~\ref{prop1-04-04-2024} and obtain the sharp decay bounds. Finally, we apply energy estimate Proposition~\ref{Energy} and obtain the refined energy estimates \eqref{Refined energy bound}.
	
For the simplicity of expression, we make the following convention:
\\
``$\lesssim$'' refers to ``smaller or equal to up to a universal constant''.
\\
``$\lesssim_{c,C_S,s_0}$" means ``smaller or equal to up to a constant determined by $(c, C_S, s_0)$''.
	\section{Preliminaries}
	\subsection{Geometry of the hyperboloidal foliation}
	We work in the (3+1)-dimensional Minkowski space-time with signature $(+, -, -, -)$ and in the Cartesian coordinates we write $(t, x) = (x^0, x^1, x^2, x^3)$, and its spacial radius is written as $r = \sqrt{(x^1)^2 + (x^2)^2 + (x^3)^2}$. We use Latin indices $a,b,c\cdots\in\{1, 2, 3\}$ for spatial components, while Greek indices $\alpha,\beta,\gamma,\cdots\in\{0, 1, 2, 3\}$ for space-time components. The indices are raised or lowered by the Minkowski metric $\eta^{\alpha\beta}$.
	
	In this article we focus on the interior of a light-cone
	$$
	\mathcal{K}: = \{(t,x)\vert r < t - 1\}.
	$$
	In the interior of this cone we recall the hyperboloidal foliation
	$$
	\mathcal{K} = \bigcup_{s > 1}^\infty\mathcal{H}_s^*,
	$$
	where
	$$
	\mathcal{H}_s := \{(t,x)\vert t^2 - r^2 = s^2, t > 0\},\quad \mathcal{H}_s^* := \mathcal{H}_s \cap \mathcal{K}.
	$$
	We also denote by $\mathcal{K}_{[s_0,s_1]} := \{(t,x) \vert s_0^2\leqslant t^2-r^2 \leqslant s_1^2,\  r<t-1\}$ the subdomain of $\Kcal$ limited by $\Hcal_{s_0}$ and $\Hcal_{s_1}$, and $\del\mathcal{K}_{[s_0,s_1]}:= \{r= t-1\vert \sqrt{r^2 + s_0^2}\leqslant t\leqslant \sqrt{r^2+s_1^2}\}$.
	
	With respect to the Euclidean metric of $\RR^{1+3}$, the normal vector and the volume form of $\Hcal_s$ is written as:
	$$
	\vec{n} = \frac{1}{\sqrt{t^2+r^2}}(t,-x^a),\quad d\sigma = \sqrt{1+(r/t)^2} dx,
	$$
	which leads to
	\begin{equation}
		\vec{n}d\sigma = (1, - x^a/t)dx.
	\end{equation}
	
	\subsection{Frames of interest}
	We recall the semi-hyperboloidal frame:
	\begin{equation}
		\underline\del_0 := \del_t,\qquad \underline\del_{a} := (x^a/t)\del_{t}+\del_{a}.
	\end{equation}
	The transition matrices between this frame and the natural frame $\del_\beta$ are, in the sense of $\delu_\alpha = \Phi^\beta_\alpha\del_\beta$ and $\del_\alpha = \Psi^\beta_\alpha\delu_\beta$,
	\begin{equation}\label{eq1-30-10-2024}
		\Phi = 
		\left(
		\begin{array}{cccc}
			1 & 0 & 0 & 0\\
			x^1/t & 1 & 0 & 0\\
			x^2/t & 0 & 1 & 0\\
			x^3/t & 0 & 0 & 1\\
		\end{array}
		\right),
		\quad\quad
		\Psi := \Phi^{-1} = 
		\left(
		\begin{array}{cccc}
			1 & 0 & 0 & 0\\
			-x^1/t & 1 & 0 & 0\\
			-x^2/t & 0 & 1 & 0\\
			-x^3/t & 0 & 0 & 1\\
		\end{array}
		\right),
	\end{equation}
	These matrices are smooth within the cone $\mathcal{K}$.
	Notice that the vector fields $\underline\del_a$ generate the tangent space of the hyperboloids, and the normal vector of hyperboloids with respect to the Minkowski metric can be written as $\underline\del_{\perp} := (t/s)\del_{t} + (x^a/s)\del_{a}$.
	
	For a quadratic form $A$, one can write it within the semi-hyperboloid frame:
	\begin{equation}
		A^{\alpha\beta}\del_\alpha\otimes\del_\beta = \underline{A}^{\alpha\beta}\delu_\alpha\otimes\delu_\beta
	\end{equation}
	where $A^{\alpha\beta} = \underline{A}^{\alpha'\beta'}\Phi_{\alpha'}^\alpha\Phi_{\beta'}^\beta$.
	In particular, $\underline{A}^{00}$ is expressed in \eqref{eq5-29-08-2024}.
	The corresponding semi-hyperboloidal co-frame can be represented as:
	\begin{equation}
		\theta^0 := dt - \sum_a(x^a/t)dx^a, \qquad \theta^{a} := dx^a.
	\end{equation}
	In the semi-hyperboloidal frame, the Minkowski metric is written as:
	\begin{equation}
		\big(\etau _{\alpha\beta}\big)=
		\left(
		\begin{array}{cccc}
			1 & x^1/t      & x^1/t & x^3/t \\
			x^1/t& (x^1/t)^2-1 & x^1x^2/t^2 & x^1x^3/t^2\\
			x^2/t& x^2x^1/t^2 & (x^2/t)^2-1 & x^2x^3/t^2 \\
			x^3/t& x^3x^1/t^2 & x^3x^2/t^2 & (x^3/t)^2-1
		\end{array}
		\right),    
	\end{equation}
	\begin{equation}
		\big(\etau ^{\alpha\beta}\big)=
		\left(
		\begin{array}{cccc}
			(s/t)^2 & x^1/t      & x^1/t & x^3/t \\
			x^1/t& -1 & 0 & 0\\
			x^2/t& 0 & -1 & 0\\
			x^3/t& 0 & 0 & -1
		\end{array}
		\right).      
	\end{equation}
	
We also recall the hyperbolic variables:
$$
s := \sqrt{t^2-r^2},\quad x^a=x^a
$$	
which is an other parameterization of $\Kcal$. We denote by $\{\delb_{\alpha}\}$ the associate natural frame: 
$$
\delb_0 = \delb_s = (s/t)\del_t,\quad \delb_a = (x^a/t)\del_t + \del_a = \delu_a.
$$

	\subsection{The energy estimates on hyperboloids}
	For any functions $u$ defined in (an open subset of) $\mathcal{K} = \{r<t-1\}$, we define their integral on the hyperboloids as:
	\begin{equation}
		\|u\|_{L^1_f(\mathcal{H}_s)} := \int_{\mathcal{H}_s} u dx =\int_{\mathbb{R}^3}u\big(\sqrt{s^2+r^2},x\big)dx.
	\end{equation}
	We recall the following standard energy:
	\begin{equation}\label{eq1-21-04-2023}
		\Ebf_{c}(s,u)=\int_{\mathcal{H}_s}\Big(|\del_t u|^2+\sum_{a}|\del_{a}u|^2+2(x^a/t)\del_tu\del_{a}u+c^2u^2 \Big)dx.
	\end{equation}
	The hyperboloidal conformal energy is defined as 
	\begin{equation}
		\begin{aligned}
			&\Ebf_{\text{con}}(s,u):= \int_{\mathcal{H}_s}\Big((Ku+2u)^2+\sum_{a}(s\delu_a u)^2 \Big)dx\\
		\end{aligned} 
	\end{equation}
	where 
	$$
	Ku=(s\del_s+2x^a\delu_a)u.
	$$
	They both satisfy the following energy estimate(cf.\cite[Proposition~2.2]{Y. Ma-2021-strongcoupling}, see also \cite{W. Wong-2017})
	\begin{proposition}\label{Energy}
		For any function $u$ defined in $\mathcal{K}_{[s_0,s_1]}$ and vanishes near $\del\mathcal{K}$, for all $s\in[s_0,s_1]$,
		\begin{equation}\label{eq1-23-08-2023}
			\Ebf_{c}(s,u)^{1/2}\leqslant \Ebf_{c}(2,u)^{1/2}+C\int_{s_0}^s\|\square u\|_{L^2_f(\mathcal{H}_{\bar s})}d\bar s,
		\end{equation}
		\begin{equation}\label{eq2-23-08-2023}
			\Ebf_{\con}(s,u)^{1/2}\leqslant \Ebf_{\con}(s_0, u)^{1/2}+C\int_{s_0}^s \bar s\|\square u\|_{L^2_f(\mathcal{H}_{\bar s})}d \bar s,
		\end{equation}
		where $\Box u = \eta^{\alpha\beta}\del_{\alpha}\del_{\beta}u$.
	\end{proposition}
	For the conformal energy, we have the following estimate (cf. \cite[Lemma~2.2]{Duan-Ma-Zhang-2023}):
	\begin{lemma}\label{01-08-2023}
		Let $u$ be a function defined in $\Kcal_{[s_0,s_1]}$ and vanishes near the conical boundary $\del\Kcal = \{r=t-1\}$. Then
		\begin{equation}\label{eq1-31-10-2024}
			\|s(s/t)^2\del_{\alpha}u\|_{L^2_f(\Hcal_s)} + \|(s/t)u\|_{L^2_f(\Hcal_s)}\lesssim\Ebf_{\con}(s,u)^{1/2}.
		\end{equation}
	\end{lemma}
	
	For the convenience of discussion, we introduce the following energy densities:
	\begin{equation}\label{eq2-14-08-2023}
		\ebf_c[u]:= \sum_{\alpha=1}^3(s/t)|\del_{\alpha}u|^2 +  \sum_{a=1}^3|\delu_a u|^2 + c^2u^2,
	\end{equation}
	\begin{equation}\label{eq3-14-08-2023}
		\ebf_{\con}[u] := \sum_{\alpha=0}^3|s(s/t)^2\del_{\alpha}u|^2 +  \sum_{a=1}^3|s\delu_a u|^2 + |(s/t)u|^2.
	\end{equation}
	From Lemma~\ref{01-08-2023}, it is clear that
	$$
	\int_{\Hcal_s}\ebf_{\con}[u]dx\lesssim\Ebf_{\con}(s,u).
	$$
	
	\subsection{High-order operators and Sobolev decay estimates}
	For $a=1,2,3$ we recall the Lorentz boosts:
	\begin{equation}
		L_a:= x^{a}\del_{t}+t\del_{a}=x_a\del_0-x_0\del_a.
	\end{equation}
	For a multi-index $I=(i_n,i_{n-1},\cdots, i_1)$, we note $\del^{I}:= \del_{i_n}\del_{i_{n-1}}\cdots \del_{i_1}$. Similarly, we have $L^J=L_{i_n}L_{i_{n-1}} \cdots L_{i_1}$.
	
	Let $Z$ be a high-order derivative composed by $\{\del_{\alpha}, L_a\}$. We denote by $\text{ord}(Z)$ the order of the operator, and $\text{rank}(Z)$ the number of boosts contained in $Z$. Given two integers $k \leqslant p$, it is convenient to introduce the notations: 
	\begin{equation}
		\aligned
		|u|_{p,k} &:= \max_{\ord(Z)\leqslant p\atop \rank(Z)\leqslant k}|Z u|,\quad &&|u|_p := \max_{0\leqslant k\leqslant p}|u|_{p,k},
		\\
		|\del u|_{p,k} &:= \max_{\alpha=0,1,2}|\del_{\alpha} u|_{p,k}, &&|\del u|_p := \max_{0\leqslant k\leqslant p}|\del u|_{p,k},
		\\
		|\del^m u|_{p,k} &:= \max_{|I|=m}|\del^I u|_{p,k}, &&|\del^m u|_p := \max_{0\leqslant k\leqslant p}|\del^I u|_{p,k},
		\\
		|\delu u|_{p,k} &:= \max_a\{|\delu_a u|_{p,k}\}, &&|\delu u|_p := \max_{0\leqslant k\leqslant p}|\delu u|_{p,k}.
		\endaligned
	\end{equation}
	We recall the following estimates based on (B.2), (2.21) and (2.22) in  \cite{Y. Ma-2021-strongcoupling}, which can be easily checked by induction:
	\begin{equation}\label{eq1-24-march-2025}
		\aligned
		|u|_{p,k}\leqslant & C\sum_{\vert I\vert +\vert J\vert \leqslant p \atop |J|\leqslant k}|\del^IL^J u|,
		\\
		|\del u|_{p,k}\leqslant & C\sum_{|I|+|J|\leqslant p\atop |J|\leqslant k,\alpha}|\del_{\alpha}\del^IL^Ju|,
		\\
		|(s/t)\del u|_{p,k}\leqslant & C(s/t)\sum_{|I|+|J|\leqslant p\atop |J|\leqslant k,\alpha}|\del_{\alpha}\del^IL^Ju|,
		\endaligned
	\end{equation} 
	\begin{equation}\label{eq2-24-march-2025}
		|\delu u|_{p,k}\leqslant C\sum_{|I| + \vert J\vert\leqslant p,a\atop |J|\leqslant k}|\delu_a\del^IL^Ju|
		+ Ct^{-1}\hspace {-0.3cm}\sum_{0\leqslant |I|+\vert J\vert \leqslant p - 1 \atop |J|\leqslant k,\alpha}\hspace {-0.3cm}|\del_{\alpha}\del^IL^Ju|\leqslant Ct^{-1}|u|_{p+1,k+1}.
	\end{equation}
	
	On the other hand, one also introduce the high-order energy densities:
	\begin{equation}
		\ebf_c^{p,k}[u] := \sum_{|I|+|J|\leqslant p\atop |J|\leqslant k}\ebf_c[\del^IL^J u],\quad
		\ebf_{\con}^{p,k}[u]:= \sum_{|I|+|J|\leqslant p\atop |J|\leqslant k}\ebf_{\con}[\del^IL^J u],
	\end{equation}
	\begin{equation}
		\ebf_c^p[u] := \sum_{k\leqslant p}\ebf_c^{p,k}[ u],\quad
		\ebf_{\con}^p[u]:= \sum_{k\leqslant p}\ebf_{\con}^{p,k}[u].
	\end{equation}
	as well as the high-order energies:
	\begin{equation}
		\Ebf_c^{p,k}(s,u) := \sum_{|I|+|J|\leqslant p\atop |J|\leqslant k}\hspace{-0.3cm} \Ebf_c(s,\del^IL^Ju),\quad
		\Ebf_{\con}^{p,k}(s,u) := \sum_{|I|+|J|\leqslant p\atop |J|\leqslant k}\hspace{-0.3cm} \Ebf_{\con}(s,\del^IL^Ju),
	\end{equation}
	\begin{equation}\label{eq5-31-12-2022}
		\Ebf_c^p(s,u) := \sum_{k\leqslant p} \Ebf_c^{p,k}(s,u),\quad
		\Ebf_{\con}^p(s,u) := \sum_{k\leqslant p} \Ebf_{\con}^{p,k}(s,u).
	\end{equation}
	
	Recall \eqref{eq2-14-08-2023} and \eqref{eq3-14-08-2023} together with \eqref{eq1-24-march-2025} and \eqref{eq2-24-march-2025},
	\begin{equation}\label{eq2-21-04-2023}
		|(s/t)\del u|_{p,k}^2 + |\delu u|_{p,k}^2 + c|u|_{p,k}^2\leqslant C\ebf_c^{p,k}[u],\quad 
	\end{equation}
	\begin{equation}\label{eq4-14-08-2023}
		|s(s/t)^2\del u|_{p,k}^2 + |s\delu u|_{p,k}^2 + |(s/t)u|_{p,k}^2\leqslant C\ebf_{\con}^{p,k}[u].
	\end{equation}

	On the other hand, we recall the following Klainerman-Sobolev type estimates established in \cite[Chapter VII]{Hörmander-book-2003}.
	\begin{proposition}\label{Klainerman-Sobolev}
		Let $u$ be a function defined in $\Kcal_{[s_0,s_1]}$ and vanishes near $\del\Kcal = \{r=t-1\}$. Then
		\begin{equation}
			\sup_{\Hcal_s}\{t^{3/2}|u|\}\leqslant C\sum_{|I|+|J|\leqslant 2}\|\del^IL^J u\|_{L^2_f(\Hcal_s)}.
		\end{equation}
	\end{proposition}
	
%
%
	Combining Proposition~\ref{Klainerman-Sobolev}
	, \eqref{eq2-21-04-2023} and \eqref{eq4-14-08-2023}, we obtain the pointwise estimates on $\Hcal_s$:
	\\
	\begin{equation}\label{eq3-21-04-2023}
		t^{3/2}(s/t)|\del u|_{p,k} + t^{3/2}|\delu u|_{p,k} + ct^{3/2}|u|_{p,k}\leqslant C\sqrt{\Ebf_c^{p+2,k+2}(s,u)},
	\end{equation}
	\begin{equation}\label{eq5-14-08-2023}
		t^{3/2}s(s/t)^2|\del u|_{p,k} + t^{3/2}s|\delu u|_{p,k} + t^{3/2}(s/t)|u|_{p,k}\leqslant C\sqrt{\Ebf_{\con}^{p+2,k+2}(s,u)}.
	\end{equation}

	\section{The $L^\infty$ estimation of wave equation and Klein-Gordon equation}
	\subsection{Estimates on wave equation}
	Given the intricate nature of system \eqref{eq01}, the Klainerman-Sobolev inequalities fail to provide adequate decay. While, we require more precise estimates that take into account the linear structure of the wave and/or Klein-Gordon equations. These estimates allow us to attain the linear decay rate even when the energies are not uniformly bounded. By leveraging these refined estimates, we can better understand the decay properties of the system and handle its complexities effectively
	\begin{proposition}[{cf. \cite{Lefloch-Ma-2016-CMP}}]\label{prop-01-10-2023}
		Let $u$ be a solution to the following Cauchy Problem:
		\begin{equation}
			\Box u=f,\quad u|_{t=2}=0,\ \partial_{t}u|_{t=2}=0,
		\end{equation}
		where the source $f$ vanishes outside of $\Kcal$, and there exists a global constant $C_f$ depending on f such that:  
		$$
		|f|\leqslant C_ft^{-2-\nu}(t-r)^{-1+\mu},\quad 0<\mu,|\nu|\leqslant 1/2.
		$$
		Then $u$ satisfies the following estimate:
		\begin{equation}
			|u(t,x)|\lesssim 
			C_f\begin{cases}
				\frac{1}{\nu\mu}(t-r)^{\mu-\nu} t^{-1}, \qquad & 0< \nu\leqslant 1/2,
				\\
				\frac{1}{|\nu|\mu}(t-r)^{\mu} t^{-1 -\nu}, &-1/2\leqslant \nu < 0.
			\end{cases}
			\label{wave component estimate}
		\end{equation}
	\end{proposition}

	\subsection{The sharp decay estimates for the Klein-Gordon equation}
	
	Consider a (sufficiently smooth and vanishing near the light cone) solution to the following linear Klein-Gordon equation
	\begin{subequations}
		\begin{align}
			&\Box v - h\del_tv + c^2v = f, \label{Klein-Gordon equation}
			\\
			&v\vert_{\mathcal{H}_2} = v_0,\quad \del_tv\vert_{\mathcal{H}_2} = v_1, \label{Klein-Gordon Cauchy}
		\end{align}
	\end{subequations}
where $h$ is a sufficiently regularity function. We remark the following calculation.
\begin{equation}\label{10-03-2023-ZWD}
\aligned
\Box v =&\, \delb_s\delb_sv + \sum_a\frac{2x^a}{s}\delb_s\delb_av - \sum_a\delb_a\delb_av + \frac{3}{s}\delb_sv
		\\
=&\, s^{-3/2}(\delb_s + (x^a/s)\delb_a)^2\big(s^{3/2}v\big) - \frac{x^ax^b}{s^2}\delb_a\delb_b v - \sum_a\delb_a\delb_a v - \frac{3x^a}{s^2}\delb_a v - \frac{3}{4s^2}v.
\endaligned
\end{equation}
Then \eqref{Klein-Gordon equation} can be written as
$$
\frac{x^ax^b}{s^2}\delb_a\delb_b v + \sum_a\delb_a\delb_a v + \frac{3x^a}{s^2}\delb_a v + \frac{3}{4s^2}v + f = s^{-3/2}\Lcal(s^{3/2}v) - h\del_t v + c^2v,
$$
where
$$
\Lcal = (\delb_s + (x^a/s)\delb_a).
$$
Thus
\begin{equation}
\aligned
\Lcal^2&(s^{3/2}v) - (t/s)h\Lcal(s^{3/2} v) + c^2s^{3/2}v 
\\
=&s^{3/2}\Big(\frac{x^ax^b}{s^2}\delb_a\delb_b v + \sum_a\delb_a\delb_a v 
+ \frac{3x^a}{s^2}\delb_a v +\frac{3}{4s^2}v\Big) 
- \frac{3t}{2s^{1/2}}hv - s^{1/2}(x^a/s)hL_a v + s^{3/2}f.
\endaligned
\end{equation}
	For convenience, we define 
	\begin{equation}\label{eq01-19-07-2024}
		w_{t,x}(\lambda):= \lambda^{3/2}v\vert_{\gamma_{t,x}(\lambda)},
	\end{equation}
	that is
	\begin{equation}\label{eq13-19-07-2024}
		\Lcal(s^{3/2}v)\vert_{\gamma_{t,x}(\lambda)} = w_{t,x}^\prime(\lambda),
	\end{equation}
	here $\gamma_{t,x} :\RR\mapsto \Kcal_{[s_0,s_1]}$ is the integral curve of $\Lcal$ with $\gamma_{t,x}(s) = (t,x), s=\sqrt{t^2-r^2}$. It is explicitly written as
	$$
	\gamma_{t,x}(\lambda) = \big\{
	\big(\lambda t/s,\lambda x/s\big)\big\}.
	$$
	Therefore we have the following result:
	\begin{lemma}\label{lemma1-22-07-2024}
		For every sufficiently smooth solution $v$ to \eqref{Klein-Gordon equation} the function \eqref{eq01-19-07-2024} satisfies the following second-order ODE
		\begin{equation}\label{eq1-10-04-2023}
			w_{t,x}''(\lambda) - D_{t,x}[h](\lambda) w_{t,x}'(\lambda) + c^2 w_{t,x}(\lambda) = R_{1,t,x}[v](\lambda) + R_{2,t,x}[h,v](\lambda)  + F_{t,x}(\lambda),
		\end{equation}
		where
		\begin{equation}\label{eq2-10-04-2023}
			\aligned
			D_{t,x}[h](\lambda) =& (t/s)h\Big|_{\gamma_{t,x}(\lambda)},
			\\
			R_{1,t,x}[v](\lambda) =& s^{3/2}\Big(\frac{x^ax^b}{s^2}\delb_a\delb_b + \sum_a\delb_a\delb_a 
			+ \frac{3x^a}{s^2}\delb_a + \frac{3}{4s^2}\Big)v \Big|_{\gamma_{t,x}(\lambda)},
			\\
			R_{2,t,x}[h,v](\lambda) =& \Big(-\frac{3t}{2s^{1/2}}hv - s^{1/2}(x^a/s)hL_a v\Big)\Big|_{\gamma_{t,x}(\lambda)},
			\\
			F_{t,x}(\lambda) = &s^{3/2}f\Big|_{\gamma_{t,x}(\lambda)}.
			\endaligned
		\end{equation}
		
		Here $\gamma_{t,x} :\RR\mapsto \Kcal_{[s_0,s_1]}$ is the integral curve of $\Lcal$. It is explicitly written as
		$$
		\gamma_{t,x}(\lambda) = \big\{
		\big(\lambda t/s,\lambda x/s\big)\big\}.
		$$
	\end{lemma}
	
%
	
	Now we need an observation on ODE.
	\begin{lemma}\label{lemma1-19-07-2024}
		Let $v$ be the solution to the following ODE on the interval $[\lambda_0,\lambda_1]$:
		\begin{equation}\label{eq1-11-03-2023}
			w''(\lambda) - D(\lambda)w'(\lambda) + c^2w(\lambda) = f(\lambda). 
		\end{equation}
		Suppose that $\sup_{\lambda\in[\lambda_0,\lambda_1]}|D(\lambda)|\leqslant c$, and $\exists S: [\lambda_0,\lambda_1]\mapsto\RR$, satisfying $S(\eta) + D(\eta)\leqslant 0$, and $\int_{\lambda_0}^{\lambda} |S(\eta)|d\eta\leqslant C_S$, where $C_S$ is a constant. Then
		\begin{equation}\label{eq12-19-07-2024}
			|w(\lambda)| + |w^\prime(\lambda)|\lesssim_{c,C_S}  |w(\lambda_0)| + |w^\prime(\lambda_0)| 
			+ \int_{\lambda_0}^\lambda \vert f(\tau)\vert d\tau.
		\end{equation}
	\end{lemma}
	\begin{proof}
	The equation \eqref{eq1-11-03-2023} can be written into a first-order linear system
		\begin{equation}\label{eq2-11-03-2023}
			\left(
			\begin{array}{c}
				w^\prime
				\\
				w
			\end{array}
			\right)'
			-
			A
			\left(
			\begin{array}{c}
				w^\prime
				\\
				w
			\end{array}
			\right)
			= 
			\left(
			\begin{array}{c}
				f
				\\
				0
			\end{array}
			\right),
		\end{equation}
		in which $A = A(\lambda)$ can be diagonalized
		\begin{equation}
			A(\lambda) = P(\lambda)\text{diag}(p_+,p_-)P^{-1}(\lambda),
		\end{equation}
		here 
		\begin{equation}
			A(\lambda) = \left(
			\begin{array}{cc}
				D & -c^2
				\\
				1 & 0
			\end{array}
			\right),
			\quad
			P(\lambda)=
			\left(
			\begin{array}{cc}
				p_+ & p_-
				\\
				1 &1
			\end{array}
			\right),
		\end{equation}
	and
	\begin{equation}
			P^{-1}(\lambda)=(p_+-p_-)^{-1}
			\left(
			\begin{array}{cc}
			1 & -p_-
			\\
			-1 &p_+
			\end{array}
			\right)
		\end{equation}
		with $ p_{+,-} = \frac{D\pm \sqrt{D^2-4c^2}}{2}$.
		
		Remark that $p_+-p_- = \sqrt{D^2-4c^2}$. Then we have, thanks to the condition $\sup_{\lambda\in[\lambda_0,\lambda_1]}|D(\lambda)|\leqslant c$,
		$$
		|P(\lambda)| + |P^{-1}(\lambda)|\lesssim_c 1.
		$$
		Now \eqref{eq2-11-03-2023} is written as
		$$
		P^{-1}(\lambda)\left(
		\begin{array}{c}
			w'
			\\
			w
		\end{array}
		\right)'
		-
		\left(
		\begin{array}{cc}
			p_+ &0
			\\
			0 &p_-
		\end{array}
		\right)
		P^{-1}(\lambda)\left(
		\begin{array}{c}
			w'
			\\
			w
		\end{array}
		\right)
		= 
		P^{-1}(\lambda)\left(
		\begin{array}{c}
			f
			\\
			0
		\end{array}
		\right).
		$$
		We set $\tilde{W}(\lambda) = P^{-1}(\lambda)\left(
		\begin{array}{c}
			w'
			\\
			w
		\end{array}
		\right)$, then
		\begin{equation}
			\tilde{W}'(\lambda) 
			-
			\left(
			\begin{array}{cc}
				p_+ &0
				\\
				0 &p_-
			\end{array}
			\right)
			\tilde{W}(\lambda)
			=
			F(\lambda),
		\end{equation}
		where
		\begin{equation}\label{eq1-22-05-2023}
			F(\lambda) = 
			P^{-1}(\lambda)\left(
			\begin{array}{c}
				f
				\\
				0
			\end{array}
			\right)
			+
			(P^{-1}(\lambda))'\left(
			\begin{array}{c}
				w'
				\\
				w
			\end{array}
			\right).
		\end{equation}
		We denote by
		\begin{equation}
			\tilde{W} = 
			\left(
			\begin{array}{cc}
				\tilde{W}_+
				\\
				\tilde{W}_-
			\end{array}
			\right),
			\quad
			F = \left(
			\begin{array}{cc}
				F_+
				\\
				F_-
			\end{array}
			\right).
		\end{equation}
		Then
		\begin{equation}
			\aligned
			\tilde{W}_\pm(\lambda) =& \tilde{W}_\pm(\lambda_0)e^{\int_{\lambda_0}^\lambda p_\pm(\eta)d\eta} + \int_{\lambda_0}^\lambda F_\pm(\tau)e^{\int_{\tau}^{\lambda}p_\pm(\eta)d\eta}d\tau.
			\endaligned
			\label{eq1-01-04-2023}
		\end{equation}
		Here we notice that
		\begin{equation}\label{eq3-10-04-2023}
			\int_{\tau}^\lambda p_{\pm}(\eta)d\eta 
			=  \frac{1}{2}\int_{\tau}^\lambda D(\eta)d\eta 
			\pm \frac{i}{2}\int_{\tau}^{\lambda} \sqrt{4c^2-D^2(\eta)}d\eta
		\end{equation}
		To achieve our purpose, we recall the barrier function $S(\eta)$ which satisfies
		\begin{equation}\label{eq5-20-04-2023}
			D(\eta) + S(\eta) \leqslant 0 \quad\forall\eta\in[\lambda_0,\lambda_1],
		\end{equation}
		and notice that
		\begin{equation}\label{eq1-05-05-2023}
			\int_{\lambda_0}^{\lambda_1}S(\eta)d\eta\leqslant C_S,
		\end{equation}
		where $C_S$ is a constant. Then \eqref{eq3-10-04-2023} is written as follow
		\begin{equation}\label{eq6-20-04-2023}
			\aligned
			\int_{\tau}^\lambda p_{\pm}(\eta)d\eta 
			= & \frac{1}{2}\int_{\tau}^\lambda D(\eta)d\eta 
			\pm \frac{i}{2}\int_{\tau}^{\lambda} \sqrt{4c^2-D^2(\eta)}d\eta\\
			= & \frac{1}{2}\int_{\tau}^\lambda \big(D(\eta) + S(\eta)\big)d\eta - \frac{1}{2}\int_{\tau}^\lambda S(\eta)d\eta \\
			& \pm \frac{i}{2}\int_{\tau}^{\lambda} \sqrt{4c^2 - D^2(\eta)}d\eta.
			\endaligned
		\end{equation}
		
		We only need to consider if $S(\lambda)$ is integrable, which we have guaranteed in \eqref{eq1-05-05-2023}. More precisely:
		\begin{equation}\label{eq1-13-05-2023}
			\aligned
			e^{\int_{\lambda_0}^\lambda p_\pm(\eta)d\eta} 
			& = e^{\frac{1}{2}\int_{\tau}^\lambda (D(\eta) - S(\eta))d\eta}
			\cdot e^{-\frac{1}{2}\int_{\tau}^\lambda S(\eta)d\eta} \cdot e^{\pm\frac{i}{2}\int_{\tau}^{\lambda} \sqrt{4c^2 - D^2(\eta)}d\eta}\\
			& \leqslant 1 \cdot e^{\frac{1}{2}C_S}\cdot e^{\pm\frac{i}{2}\int_{\tau}^{\lambda} \sqrt{4c^2-D^2(\eta)}d\eta}\\
			&\leqslant e^{\frac{1}{2}C_S\pm\frac{i}{2}\int_{\tau}^{\lambda} \sqrt{4c^2- D^2(\eta)}d\eta}
			 \lesssim_{c,C_S} 1.
			\endaligned
		\end{equation}
		
		Notice that $P(\lambda)$ and its inverse matrix $\big(P(\lambda)\big)^{-1}$ are uniformly bounded. Combine \eqref{eq1-22-05-2023}, \eqref{eq1-01-04-2023} and \eqref{eq1-13-05-2023} we achieve our target.
	\end{proof}
	
	\begin{proposition}\label{prop1-04-04-2024}
		Considering a sufficiently smooth solution $v$ to the Klein-Gordon problem \eqref{Klein-Gordon equation} in $\mathcal{K}_{[s_0,s_1]}$ and vanishing near the light cone. Suppose that for all $(t,x)\in\mathcal{K}_{[s_0,s_1]}$, one has $\vert (t/s)h\vert\leqslant c$, and there exists a function $S$, s.t. $\forall(t,x)\in\Kcal_{[s_0,s_1]}$, $S(t,x)+(t/s)h\leqslant 0$, and 
		$$
		\int_{\lambda_0}^{s_1}\big|S\vert_{\gamma_{t,x}}(\lambda)\big|d\lambda\leqslant C_S
		$$ 
		with $C_S$ a constant. Then for any $\eta \in\RR$, the following estimate holds:
		\begin{equation}\label{eq1-04-04-2024}
			(s/t)^{\eta} s^{3/2}\big(|v(t,x)| + (s/t)|\del v(t,x)| \big)\lesssim_{c,C_S,s_0} V(t,x),
		\end{equation}
		where the function $V$ is defined by distinguishing between the region ``near" and ``far" from the light cone:
		\begin{equation}\label{Function V}
			V(t,x):= 
			\left\{
			\aligned
			&
			(s/t)^{\eta}\sup_{\mathcal{H}_{s_0}}{(|v(t,x)| + |\del v(t,x)|)} + (s/t)^{\eta}s^{1/2}\vert v\vert_1(t,x)
			\\
			&\quad + (s/t)^{\eta}\int_{\lambda_0}^s \vert R[h,v] + \tau^{3/2} f\vert\Big\vert_{\gamma_{t,x}(\tau)} d\tau, 
			&&
			0\leqslant r/t\leqslant \frac{s_0^2-1}{s_0^2+1},
			\\
			&
			(s/t)^{\eta}s^{1/2}\vert v\vert_1(t,x) + (s/t)^{\eta}\int_{\lambda_0}^s \vert R[h,v] + \tau^{3/2} f\vert\Big\vert_{\gamma_{t,x}(\tau)} d\tau, 
			&&\frac{s_0^2-1}{s_0^2+1}\leqslant r/t < 1,
			\endaligned
			\right.
		\end{equation}
		with
		\begin{equation}
			\lambda_0 =
			\begin{cases}
				s_0,\quad &0\leqslant r/t\leqslant \frac{s_0^2-1}{s_0^2+1},\\
				\sqrt{\frac{t+r}{t-r}}\sim (t/s),\quad &\frac{s_0^2-1}{s_0^2+1}\leqslant r/t\leqslant 1,
			\end{cases}
		\end{equation}
		and
		\begin{equation}
		\aligned
			R[h,v] =& s^{3/2}\Big(\frac{x^ax^b}{s^2}\delb_a\delb_b v + \sum_a\delb_a\delb_a v 
			+ \frac{3x^a}{s^2}\delb_a v + \frac{3}{4s^2}v\Big) 
			\\
			&+  \Big(-\frac{3}{2}s^{1/2}hv - s^{1/2}(x^a/s)hL_a v\Big).
		\endaligned
		\end{equation}
		Here ``$\lesssim_{c,C_S,s_0}$" means smaller or equal to up to a constant determined by $(c, C_S, s_0)$.
	\end{proposition}	

\begin{remark}
In \eqref{eq1-04-04-2024} we multiplied the estimate by a factor $(s/t)^{\eta}$. This is because in the region where $0\leqslant r/t\leqslant \frac{s_0^2-1}{s_0^2+1}$, one has $\frac{2s_0}{s_0^2+1}\leqslant (s/t)\leqslant 1$ which is uniformly bounded from infinity as well as from zero. Thus $(s/t)^{\eta}$ is bounded by a constant determined by $\eta$ for both positive and negative $\eta$ (we prefer in fact the negative ones). In the region where $\frac{s_0^2-1}{s_0^2+1}\leqslant r/t < 1$, one no longer has $(s/t)$ being uniformly far from zero, thus $(s/t)^{\eta}$ is no longer bounded for negative $\eta$. However, due to the fact that the interval of integration in  \eqref{Function V} for this region starts at $\lambda = \lambda_0 \sim (s/t)^{-1}$, this will compensate the singularity of $(s/t)^{\eta}$ (for $|\eta|$ not too large).
\end{remark}

	\begin{proof}
		Recalling the definition of $w_{t,x}(\lambda)$ in \eqref{eq01-19-07-2024}, when $\lambda\geq s_0\geq 2$ we have
		\begin{equation}\label{eq19-19-07-2024}
			 \lambda^{3/2}\big(\vert v(\lambda t/s,\lambda x/s) + \vert\Lcal v(\lambda t/s,\lambda x/s)\vert\big)\lesssim\vert w^\prime_{t,x}(\lambda)\vert + \vert w_{t,x}(\lambda)\vert,
		\end{equation}
		where $\Lcal = \delb_s + (x^a/s)\delb_a$. We also observe that
		\begin{equation}\label{eq20-19-07-2024}
			\Lcal = (s/t)\del_t + s^{-1}(x^a/t)L_a,
		\end{equation}
		thus
		\begin{equation}\label{eq21-19-07-2024}
			\aligned
			\lambda^{3/2}&\big(\vert v(\lambda t/s,\lambda x/s)\vert + (s/t)\vert\del v(\lambda t/s,\lambda x/s)\vert\big)
			\\
			&\lesssim \lambda^{-1}\vert v\vert_1(\lambda t/s,\lambda x/s) + \vert w_{t,x}^\prime(\lambda)\vert + \vert w_{t,x}(\lambda)\vert.
			\endaligned
		\end{equation}
		By combining \eqref{eq1-10-04-2023} and \eqref{eq12-19-07-2024}, we derive that
		\begin{equation}
		\aligned
		|w_{t,x}(\lambda)| + |w_{t,x}^\prime(\lambda)|
		\lesssim_{c,C_S}&  |w_{t,x}(\lambda_0)| + |w_{t,x}^\prime(\lambda_0)| 
		\\	
		&+ \int_{\lambda_0}^\lambda \vert R[h,v] + \tau ^{3/2}f\vert \Big\vert_{\gamma_{t,x}(\tau)}d\tau.
		\endaligned
		\end{equation}
		Recall that $v$ vanishes near $\del\mathcal{K}$, then when $\frac{s_0^2-1}{s_0^2+1}\leqslant r/t\leqslant 1$, $w_{t,x}(\lambda_0) = w_{t,x}^\prime(\lambda_0)= 0$. When $0\leqslant r/t\leqslant \frac{s_0^2-1}{s_0^2+1}$, $w_{t,x}(\lambda_0)$ is determined by the restriction of $v$ and $\del v$ on $\mathcal{H}_{s_0}$. This gives the desired result.

	\end{proof}

	\subsection{Decomposition of the wave component}\label{subsection3.3}
	In order to apply Proposition~\ref{prop1-04-04-2024}, we need to construct the barrier function $S(t,x)$. To this purpose, we are going to establish the following decomposition:
	\begin{equation}
		u+\frac{B}{2c^2}v^2 = \overbrace{u_{L} + u_g}^{\text{sufficient decay}} + \underbrace{u_b + \frac{B}{2c^2}v^2}_{\text{good sign}} .
	\end{equation}

	To construct this decomposition, we first remark that
	\begin{equation}\label{eq1.6}
		u = u_L + \tilde{u},
	\end{equation}
	where:
	\begin{equation}\label{eq1-22-04-2023}
		\begin{cases}
			&\Box u_L = 0,
			\\
			&u_L\vert_{t=2} = u_0+\kappa v_0^2,
			\\
			& \del_tu_L\vert_{t=2} = u_1 + 2\kappa v_0v_1,
		\end{cases}
	\end{equation}
  with $\kappa$ is a pending constant. That is, $u_L$ is the initial-data-contribution, and
	\begin{equation}\label{eq2-22-04-2023}
		\begin{cases}
			& \Box \tilde{u} = A^{\alpha\beta}\del_\alpha v\del_\beta v + Bv^2,
			\\
			& \tilde{u}\vert_{t=2} = -\kappa v_0^2,
			\\
			& \del_t\tilde{u}\vert_{t=2} = -2\kappa v_0v_1.
		\end{cases}
	\end{equation}
	We first remark that $u_L$ is the solution of the homogeneous equation, naturally it has enough decay.
	
	Consider $w = \tilde{u} + \kappa v^2$, then 
	\begin{equation}\label{eq3-22-04-2023}
		\aligned
		\Box w =& \Box \tilde{u} + 2\kappa v\Box v  + 2\kappa \eta^{\alpha\beta}\del_\alpha v\del_\beta v
		\\
		=& A^{\alpha\beta}\del_{\alpha}v\del_{\beta}v + Bv^2 - 2\kappa c^2 v^2 + 2\kappa uv\del_t v + 2\kappa \eta^{\alpha\beta}\del_\alpha v\del_\beta v
		\\
		=& (\underline{A}^{00}+ 2(s/t)^2\kappa)\del_tv\del_tv  + (B - 2\kappa c^2)v^2\\
		& + \sum_{(\alpha,\beta)\neq(0,0)}(\underline{A}^{\alpha\beta} + 2\kappa \underline{\eta}^{\alpha\beta})\delu_\alpha v\delu_\beta v + 2\kappa uv\del_tv
		\endaligned
	\end{equation}
	where $\eta^{\alpha\beta}$ is the Minkowski metric, and $\underline{\eta}^{\alpha\beta}$ is the Minkowski metric in the semi-hyperboloidal frame. 
	
	To eliminate the $v^2$ term in \eqref{eq3-22-04-2023}, we chose
	\begin{equation}\label{eq2-24-04-2023}
		\kappa = \frac{B}{2c^2}.
	\end{equation}
	That is
	\begin{equation}\label{eq7-24-04-2023}
		\left\{
		\aligned
		&\Box w = \Big(\Au^{00} + \frac{B}{c^2}(s/t)^2\Big)\del_tv\del_tv  
		\\
		&\quad\quad\quad\quad+ \sum_{(\alpha,\beta)\neq(0,0)}\Big(\underline{A}^{\alpha\beta} + \frac{B}{c^2} \underline{\eta}^{\alpha\beta}\Big)\delu_\alpha v\delu_\beta v
		+ \frac{B}{c^2} uv\del_tv,
		\\
		& w\vert_{\mathcal{H}_2} = \del w\vert_{\mathcal{H}_2} = 0.
		\endaligned
		\right.
	\end{equation}
	Then let us split $w$ in two parts
	\begin{equation}
		w = u_g + \Big(u_b + \frac{B}{2c^2}v^2\Big),
	\end{equation}
	where
	\begin{equation}\label{eq4-22-04-2023}
\left\{
\aligned
			&\Box u_g = \sum_{(\alpha,\beta)\neq(0,0)}(\underline{A}^{\alpha\beta} + \frac{B}{c^2} \underline{\eta}^{\alpha\beta})\delu_\alpha v\delu_\beta v + \frac{B}{c^2} uv\del_tv,
			\\
			& u_g\vert_{\mathcal{H}_2} = \del_tu_g\vert_{\mathcal{H}_2} = 0,
\endaligned
\right.
	\end{equation}
	and
	\begin{equation}\label{eq4-22-07-2024}
\left\{
\aligned
			&\Box (u_b + \frac{B}{2c^2}v^2) = (\underline{A}^{00} + \frac{B}{c^2}(s/t)^2)\del_tv\del_tv, 
			\\
			& (u_b + \frac{B}{2c^2} v^2)\big\vert_{\mathcal{H}_2} = \del_t(u_b + \frac{B}{2c^2} v^2)\big\vert_{\mathcal{H}_2} = 0.
\endaligned
\right.
	\end{equation}
	The quadratic source terms in \eqref{eq4-22-04-2023} content at least one hyperbolic derivative, thus $u_g$ enjoy sufficient decay.
	
	For $u_b$, instead of decay rate, we demand that it is of good sign. Recall that in dimension 3+1 the fundamental solution of the wave equation is a positive measure. Thus \eqref{eq1-24-07-2024} guarantees 
\begin{equation}\label{eq4-20-feb-2025}
u_b + \frac{B}{2c^2}v^2\leqslant 0. 
\end{equation}

Now we make the following conclusion:
\begin{equation}
u = u_L + u_g + u_b
\end{equation}
where $u_L + u_g$ enjoy sufficient decay (to be precised in Section~\ref{subsec1-feb-20-2025}), while $u_b$ can be controlled by the barrier function 
\begin{equation}
S(t,x) := \frac{B}{2c^2}v^2
\end{equation}
in the sens of \eqref{eq4-20-feb-2025}.

	\section{Global existence: direct estimate}
	\subsection{The bootstrap argument}
	
	From this section onwards, we begin the proof of Theorem \ref{theorem 1.1}, which is a bootstrap argument along the lines of the method presented in \cite{Lefloch-Ma-2014}. We first summarize our strategy on following two observations:
	\\
	1.The local solution to equation \eqref{eq01} cannot approach its maximal time of existence, denoted as $s^*$, with bounded energy of sufficiently high-order. If such an approach were possible, we could apply the local existence theory and construct a local solution to equation \eqref{eq01} starting from $(s^*-\eps)$ with initial data equal to the restriction of the local solution at time $(s^*-\eps)$. This construction would enable us to extend the local solution to $(s^*-\eps+\delta)$, where $\delta$ is determined by the system itself and the high-order energy bounds (independent of $\eps$). However, when $\eps<\delta$, the local solution would be extended beyond $s^*$, contradicting the fact that $s^*$ is the maximal time of existence.
	\\
	2.The high-order energies remain continuous with respect to the time variable whenever the local solution exists. This continuity is a direct consequence of the local existence theory.
	
	Based on the above observations, and suppose that the local solution $(u,v)$ to \eqref{eq01} satisfies a set of high-order energy bounds on time interval $[s_0,s_1]$ (contained in the maximal interval of existence). If we can show that the same set of energies satisfies strictly stronger bounds on the same interval, then one concludes that this local solution extends to time infinity. To see this, suppose that\footnote{The existence and smallness of the solution in the region limited by $\del\Kcal, \{t=2\}$ and $\Hcal_2$, i.e., $\{(t,x)\in\Kcal|2\leqslant t\leqslant \sqrt{4-r^2}\}$ can be easily obtained by local theory.}
	\begin{equation}\label{eq2-20-07-2023}
		\aligned
		&\Ebf_{\con}^N(2,u)^{1/2}\leqslant C_0\varepsilon,\\
		&\Ebf_c^N(2,v)^{1/2}\leqslant C_0\varepsilon.
		\endaligned
	\end{equation}
	Let $[2,s_1]$ be the maximal time interval in which the following energy estimate holds:
	\begin{equation}
		\begin{aligned}
			&\Ebf_{\con}^N(s,u)^{1/2}\leqslant C_1\varepsilon s^{1/2 + \delta},
			\\
			&\Ebf_c^N(s,v)^{1/2}\leqslant C_1\varepsilon s^\delta
			\label{Energy bound}
		\end{aligned}
	\end{equation}
	where $C_1>C_0$ is sufficiently large, and $\delta>0$ be determined later. Then by continuity, when $s=s_1$, at least one of \eqref{Energy bound} becomes equality. However, if we can show that (based on \eqref{Energy bound}) 
	\begin{equation}
		\begin{aligned}
			&\Ebf_{\con}^N (s, u)^{1/2}\leqslant (1/2)C_1\varepsilon s^{1/2 + \delta},
			\\
			&\Ebf_c^N(s, v)^{1/2}\leqslant (1/2)C_1\varepsilon s^\delta,
			\label{Refined energy bound}
		\end{aligned}
	\end{equation}
	then we conclude that \eqref{Energy bound} holds on $[2, s^* )$. However this is impossible when $s^*< \infty$ due to the first observation. We thus obtain the desired global-in-time existence. Therefore we need to establish the following result:
	\begin{proposition}\label{Bootstrap theorem}
		Let $N\geqslant 7$ and $0<\delta<1/10$. Suppose that \eqref{Energy bound} holds on $[2,s_1]$. Then for $C_1>2C_0$ and $\vep$ sufficiently small, \eqref{Refined energy bound} holds on the same time interval.
	\end{proposition}
	
	The rest of this article is mainly devoted to the proof of the above Proposition. In the following discussion, we apply the expression $A\lesssim B$ for a inequality $A\leqslant CB$ with $C$ a constant determined by $\delta, N$ and the system \eqref{eq01}.
	
	\subsection{Direct $L^2$ estimate and pointwise estimates}
	Based on \eqref{eq2-21-04-2023}, \eqref{eq4-14-08-2023} and \eqref{Energy bound}, we have the following $L^2$ estimates:
	\begin{equation}\label{eq1-20-08-2024}
		\| (s/t) |u|_N\|_{L^2_f(\mathcal{H}_s)} + \| s |\delu u|_N\|_{L^2_f(\mathcal{H}_s)} + \|s(s/t)^2|\del u|_N\|_{L^2_f(\mathcal{H}_s)}\lesssim C_1 \varepsilon s^{1/2+\delta},
	\end{equation}
	\begin{equation}\label{eq2-20-08-2024}
		\|c|v|_N\|_{L^2_f(\Hcal_s)} + \|(s/t)|\del v|_N\|_{L^2_f(\Hcal_s)} + \||\delu v|_N\|_{L^2_f(\Hcal_s)}
		\lesssim C_1\vep s^{\delta}.
	\end{equation}
	Recalling \eqref{eq3-21-04-2023}, \eqref{eq5-14-08-2023}, \eqref{Energy bound}, one has
	\begin{equation}\label{eq3-20-08-2024}
		(s/t)|u|_{N-2} + s|\delu u|_{N-2} + s(s/t)^2|\del u|_{N-2}\lesssim C_1\vep t^{-3/2}s^{1/2+\delta},
	\end{equation}
	\begin{equation}\label{eq5-04-04-2024}
		c|v|_{N-2} + (s/t)|\del v|_{N-2} + |\delu v|_{N-2} \lesssim C_1\vep t^{-3/2}s^{\delta}.
	\end{equation}
	
	\subsection{The uniform standard energy estimations on wave component}
	We recall Proposition \ref{Energy}. For the wave equation, we apply \eqref{eq1-23-08-2023} with $c=0$:
	\begin{equation}\label{eq3-23-08-2023}
		\Ebf^N(s,u)^{1/2}\leqslant \Ebf^N(2,u)^{1/2}  + C\int^s_{2}\||\Box u|_N\|_{L^2_{f}(\mathcal{H}_\tau)}d\tau.
	\end{equation}
	Then we remark that, provided that $N\geq 5$,
	$$
	\||\del_{\alpha}v\del_{\beta}v|_{N}\|_{L^2_f(\Hcal_s)}\lesssim C_1\vep s^{\delta}\|t^{-3/2}|\del v|_N\|_{L^2_f(\Hcal_s)}\lesssim (C_1\eps)^2s^{-3/2+2\delta},
	$$
	$$
	\||v^2|_N\|_{L^2_f(\Hcal_s)}\lesssim C_1\vep s^{\delta}\|t^{-3/2}|v|_N\|_{L^2_f(\Hcal_s)}\lesssim (C_1\eps)^2s^{-3/2+2\delta}.
	$$
	Substitute these bounds into \eqref{eq3-23-08-2023}, we obtain, provided that $C_0\leqslant C_1/2$ and $C_1\vep$ sufficiently small,
	\begin{equation}\label{eq4-23-08-2023}
		\Ebf^N(s,u)^{1/2}\leqslant C_0\vep + C(C_1\vep)^2\leqslant C_1\vep,
	\end{equation}
	then we have 
	\begin{equation}\label{eq1-14-07-2025}
		\Vert (s/t)\vert\del u\vert_N\Vert_{L^2_f(\Hcal_s)} + \Vert \vert\delu_a u\Vert_{L^2_f(\Hcal_s)}\lesssim C_1\vep.
	\end{equation}
	Also this uniform energy bound, combined with \eqref{eq5-14-08-2023}, leads us to the following pointwise estimates:
	\begin{equation}\label{eq3-04-04-2024}
		(s/t)|\del_\alpha u|_{N-2} + |\delu u|_{N-2}\lesssim C_1\vep t^{-3/2}.
	\end{equation}

		\subsection{Decomposition of wave component and corresponding estimates}\label{subsec1-feb-20-2025}
		Recall the decomposition of $u$ established in Subsection \ref{subsection3.3}
		\begin{equation}\label{eq3-20-feb-2025}
			u = u_{L} + u_g + u_b .
		\end{equation}
		In this subsection, we will establish the following estimates:
		\begin{equation}\label{eq9-24-07-2024}
			\aligned
			&\vert u_L\vert \lesssim C_0\varepsilon (s/t)^{1/2}s^{-3/2},
			\\
			&\vert u_g\vert \lesssim C_1\varepsilon (s/t)^{1/2}s^{-3/2}.
			\endaligned
		\end{equation}
\\
		We first establish the estimate on $u_g$. Remark that by \eqref{eq4-22-04-2023}, $\del^IL^Ju_g$ satisfies
		\begin{equation}\label{eq7-24-07-2024}
			\Box \del^IL^Ju_g = \del^IL^J\Big(\sum_{(\alpha,\beta)\neq(0,0)}(\underline{A}^{\alpha\beta} + \frac{B}{c^2} \underline{\eta}^{\alpha\beta})\delu_\alpha v\delu_\beta v + \frac{B}{c^2} uv\del_tv\Big),
		\end{equation}
	
		Our idea is establishing a uniform conformal energy estimate of  \eqref{eq7-24-07-2024}. To do so, we need to calculate the $L^2$ bound of the right-hand side of \eqref{eq7-24-07-2024}.		
		We have following estimates:
		\begin{lemma}
			Following the above notation, for $I+J\leqslant N-2$
			\begin{equation}\label{eq2-20-feb-2025}
				s\Vert \Box \del^IL^Ju_g\Vert_{L_f^2(\mathcal{H}_s)} \lesssim (C_1\varepsilon)^2s^{-3/2 + 3\delta}.
			\end{equation}
		\end{lemma}
		\begin{proof}
			Remark \eqref{eq1-20-08-2024} and \eqref{eq2-20-08-2024}, provided that $N\geqslant 5$, we have
			\begin{equation}\label{eq5-20-08-2024}
				\Vert\vert\delu_\alpha v\delu_\beta v\vert_{N - 1}\Vert_{L_f^2(\mathcal{H}_s)} 
				\lesssim (C_1\eps)^2s^{-5/2 + 2\delta},\quad (\alpha,\beta)\neq (0,0).
			\end{equation}
\\
To see this, we only need to recall that $\delu_av = t^{-1}L_av$. Then 
$$
|\delu_{\alpha}v\delu_{\beta}v|\lesssim \sum_{1\leqslant p\leqslant N-1}|t^{-1}L v|_{N-p}|\del v|_{p}.
$$
We remark that (due to the homogeneity of $L_a$ and $\del_{\alpha}$)
$$
|t^{-1}Lv|_{N-p}\lesssim t^{-1}|Lv|_{N-p}\lesssim t^{-1}|v|_{N+1-p}.
$$
Then, also by \eqref{eq2-20-08-2024} and \eqref{eq5-04-04-2024},
$$
\aligned
\||t^{-1}Lv|_{N-p}|\del v|_p\|_{L_f^2(\Hcal_s)}\lesssim&
\begin{cases}
\|t^{-1}|v|_{N}|\del v|_{N-2}\|_{L_f^2(\Hcal_s)},\quad &1\leqslant p\leqslant N-2,
\\
\|t^{-1}|v|_{N-2}|\del v|_{N-1}\|_{L_f^2(\Hcal_s)},\quad & p= N - 1 .
\end{cases}
\\
\lesssim&
\begin{cases}
C_1\eps s^{-5/2+\delta}\||v|_N\|_{L_f^2(\Hcal_s)}
\\
C_1\eps s^{-5/2+\delta}\|(s/t)|\del v|_N\|_{L_f^2(\Hcal_s)}
\end{cases}
\\
\lesssim &
\begin{cases}
(C_1\eps)^2s^{-5/2+2\delta}
\\
(C_1\eps)^2s^{-5/2 + 2\delta}.
\end{cases}
\endaligned
$$
Similarly, 
\begin{equation}\label{eq4-20-08-2024}
\aligned
\Vert \vert uv\del_tv\vert_{N - 2}\Vert_{L_f^2(\mathcal{H}_s)} 
&\lesssim C_1\varepsilon s^{1/2+\delta}\Vert (t/s)t^{-3/2}v\del_tv \Vert_{L_f^2(\mathcal{H}_s)} 
\\
& \lesssim (C_1\varepsilon)^2s^{1/2+2\delta} \Vert (t/s) t^{-3}\del_tv\Vert_{L_f^2(\mathcal{H}_s)} \lesssim (C_1\varepsilon)^3s^{-5/2+3\delta}.
\endaligned
\end{equation}
Combine \eqref{eq5-20-08-2024} and \eqref{eq4-20-08-2024}, we derive \eqref{eq2-20-feb-2025}.
\end{proof}
Form \eqref{eq2-20-feb-2025}, we easily get the following bound via Proposition~\ref{Energy}, estimate \eqref{eq2-23-08-2023}:
		\begin{proposition}
			For \eqref{eq7-24-07-2024} we have following uniform energy estimate:
			\begin{equation}
				\Ebf_\con^{N-2}(s,u_g)\lesssim C_1\varepsilon.
			\end{equation}
		\end{proposition}
		
		\begin{corollary}\label{cor1-24-07-2024}
			Combine this uniform energy bound with \eqref{eq5-14-08-2023}, leads us to the following pointwise estimates:
			\begin{equation}\label{eq6-20-feb-2025}
				(s/t)|u_g|_{N-4} + s|\delu u_g|_{N-4} + s(s/t)^2|\del u_g|_{N-4}\lesssim C_1\vep t^{-3/2}.
			\end{equation}
		\end{corollary}
		
		For $u_L$, we remark that it is the solution of a free-linear wave equation with compactly supported initial data on initial slits. Which means the conformal energy estimate of $u_L$ of every order are conserved. Finally, we arrive \eqref{eq9-24-07-2024}.

	\section{The sharp decay estimates}
	In this section, we establish the following estimates,
	\begin{equation}\label{eq11-21-08-2024}
		\vert u\vert_k\lesssim
		\left\{
		\begin{array}{cc}
			\begin{aligned}
				&C_1\varepsilon t^{-1},  &&k=0,
				\\
				&C_1\varepsilon t^{-1}s^{CC_1\varepsilon}, &&1\leqslant k\leqslant N-4,
			\end{aligned}
		\end{array}
		\right.
	\end{equation}
	\begin{equation}\label{eq12-21-08-2024}
		\vert v\vert _{N-4,k} + (s/t)\vert \del v\vert _{N-4,k}\lesssim
		\left\{
		\begin{array}{cc}
			\begin{aligned}
				&C_1\varepsilon(s/t)^{2-2\delta}s^{-3/2}, &&k=0,
				\\
				&C_1\varepsilon(s/t)^{2-2\delta}s^{-3/2+CC_1\varepsilon},&&1\leqslant k\leqslant N-4,
			\end{aligned}
		\end{array}
		\right.
	\end{equation}
	which is usually called the sharp decay.
	
	\subsection{The sharp decay estimates for the Klein-Gordon equation}
	
		For the convenience of discussion, we introduce the following notations for $k\leqslant N -4$:
		\begin{equation}\label{eq1-15-04-2024}
			\aligned
			\Abf_k(s) :=& \sup_{\Kcal_{[2,s]}}\{(s/t)^{2\delta - 2}s^{3/2}(\vert v\vert_{N-4,k} + (s/t)\vert\del v\vert_{N-4,k})\},
			\\
			\Bbf_k(s) :=& \sup_{\Kcal_{[2,s]}}\{t\vert u\vert_{k,k}\}.
			\endaligned
		\end{equation}
		\\
		We are going to establish the estimate on $\Abf_k$ and $\Bbf_k$. For the Klein-Gordon component, we differentiate the Klein-Gordon equation in \eqref{eq01} with respect to $\del^IL^J$, $|I|+|J|\leqslant N-4, |J|\leqslant k$ and obtain 
		\begin{equation}\label{eq1-06-04-2024}
			\Box \del^IL^Jv + c^2\del^IL^Jv  =  \del^IL^J(u\del_tv)
		\end{equation}
Then we recall the decomposition \eqref{eq3-20-feb-2025}. This leads to
		\begin{equation}\label{eq5-24-07-2024}
			\Box (\del^IL^Jv) - u_b\del_t \del^IL^J v+ c^2(\del^IL^Jv) = [\del^IL^J,u\del_t] v + (u_g+ u_L)\del_t \del^IL^J v.
		\end{equation}
		Then we apply Proposition~\ref{prop1-04-04-2024}. To do so, we need verify the properties of the barrier function $S(t,x)$:
		\begin{equation}\label{eq2-24-07-2024}
			(t/s)u_b + \underbrace{(t/s)\frac{B}{2c^2}v^2}_{S}\leqslant 0,
		\end{equation}
		and 
		\begin{equation}\label{eq3-24-07-2024}
			\int_{\lambda_0}^{+\infty} (t/s)v^2\vert_{\gamma_{t,x}}(\lambda)d\lambda\lesssim (C_1\varepsilon)^2.
		\end{equation}
		Here \eqref{eq2-24-07-2024} is a direct result of  \eqref{eq4-20-feb-2025} as been discussed at the end of Subsection \ref{subsection3.3}. For \eqref{eq3-24-07-2024}, thank to \eqref{eq5-04-04-2024}, we have
		\begin{equation}
			\vert (t/s)v^2\vert \lesssim (C_1\varepsilon)^2t^{-2}s^{-1+2\delta},
		\end{equation}
		\begin{equation}
			\vert (t/s)v^2\vert_{\gamma_{t,x}}(\lambda)\vert \lesssim( C_1\varepsilon)^2(s/t)^2\lambda^{-3+2\delta}
		\end{equation}
		guarantee \eqref{eq3-24-07-2024}.
		
		Then we decompose
		\begin{equation}
			\aligned
			R[u_b,\del^IL^Jv] = R_1[\del^IL^Jv] + R_2[u_b,\del^IL^Jv],
			\endaligned
		\end{equation}
		where
		\begin{equation}
			R_1[\del^IL^Jv] = s^{3/2}\Big(\frac{x^ax^b}{s^2}\delb_a\delb_b  + \sum_a\delb_a\delb_a  
			+ \frac{3x^a}{s^2}\delb_a  + \frac{3}{4s^2}\Big)\del^IL^J v,
		\end{equation}
		\begin{equation}
			\aligned
			R_2[u_b,\del^IL^Jv] = &-\frac{3t}{2s^{1/2}}u_b\del^IL^Jv - s^{1/2}(x^a/s)u_bL_a \del^IL^Jv.
			\endaligned
		\end{equation}
		
		\begin{lemma}\label{Estimate of R1 and R2}
			For $\vert I\vert + \vert J\vert \leqslant N-4$, we have the following estimates in $\mathcal{K}_{[s_0,s_1]}$:
			\begin{equation}
				\aligned
				&\big\vert R_1[\del^IL^Jv] \big\vert\lesssim C_1\varepsilon(s/t)^{3/2}s^{-2+ \delta},
				\\
				&\big\vert R_2[u_b,\del^IL^Jv]\big\vert \lesssim (C_1\varepsilon)^2(s/t)s^{-2+2\delta}
				.
				\endaligned
			\end{equation}
		\end{lemma}
		\begin{proof}
			The first result has been proved in the Section 6 of  \cite{Lefloch-Ma-2016-CMP}. The essence is that the terms in $R_1$ only contain hyperbolic derivatives or enjoy decreasing coefficients. Thus we only need to bound $R_2$. To this purpose we need the following estimates:
			\begin{equation}\label{eq5-20-feb-2025}
			|u_b|_{N-2}\lesssim C_1\eps (s/t)^{1/2}s^{-1+\delta}.
			\end{equation}
			This is due to \eqref{eq3-20-08-2024}  and \eqref{eq9-24-07-2024}.
			Based on the above estimate and \eqref{eq5-04-04-2024}, we have
			\begin{equation}\label{eq12-20-08-2024}
				\aligned
				\vert ts^{-1/2}u_b\del^IL^Jv\vert &\lesssim ts^{-1/2}C_1\varepsilon (s/t)^{1/2}s^{-1+\delta}C_1\varepsilon t^{-3/2}s^\delta 
				\\
				&\lesssim (C_1\varepsilon)^2(s/t)s^{-2+2\delta},
				\endaligned
			\end{equation}
			\begin{equation}\label{eq13-20-08-2024}
				\aligned
				\vert s^{1/2}(x^a/s)u_bL_a\del^IL^Jv\vert&\lesssim s^{1/2}(t/s)C_1\varepsilon (s/t)^{1/2}s^{-1+\delta}C_1\varepsilon t^{-3/2}s^{\delta}
				\\
				&\lesssim (C_1\varepsilon)^2(s/t)s^{-2+2\delta}.
				\endaligned
			\end{equation}
			We thus obtain the estimate on $\big\vert R_2[u_b,\del^IL^Jv]\big\vert$.
		\end{proof}
		
		Then we want estimate the source terms, which are the right-hand side of \eqref{eq5-24-07-2024}. For the first term, we have 
		\begin{lemma}\label{lemma1-04-04-2024}
			Let $\vert I\vert +\vert J\vert \leqslant N - 4$ and $\vert J\vert = k$, then
			\begin{equation}\label{eq2-04-04-2024}
				\aligned
				\vert [\del^IL^J,u\del_t]v\vert\lesssim &
				 (C_1\varepsilon)^2(s/t)^2s^{-3+\delta}
				\\
				&+ (s/t)^{2-2\delta}s^{-5/2}\big(\Bbf_0(s)\Abf_{k-1}(s) + \sum_{k_1=1}^{k}\Bbf_{k_1}(s)\Abf_{k-k_1}(s)\big).
				\endaligned
			\end{equation}
			Especially when $k = 0$, the last term does not exist.
		\end{lemma}
		\begin{proof}\label{proof1-04-04-2024}
			We focus on
			\begin{equation}\label{eq4-04-04-2024}
				\aligned
				\vert[\del^IL^J,u\del_t]v\vert \lesssim 
				&\sum_{{\vert I_1\vert + \vert I_2\vert = \vert I\vert, \vert I_1\vert\geq1,} \atop\vert J_1\vert +\vert J_2\vert = \vert J\vert}\vert\del^{I_1}L^{J_1}u\vert\vert\del^{I_2}L^{J_2}\del_t v\vert\\
				& + \sum_{{\vert J_1\vert + \vert J_2\vert = \vert J\vert,}\atop \vert J_1\vert\geq1}\vert L^{J_1}u\vert\vert\del^{I}L^{J_2}\del_tv\vert + \vert u\vert\Big\vert[\del^IL^J,\del_t]v\Big\vert.
				\endaligned
				\end{equation}
		Here we notice that the last two terms do not exist when $k = 0$. For the first term, we remark that, thanks to \eqref{eq5-04-04-2024} and \eqref{eq3-04-04-2024},
		\begin{equation}\label{eq6-04-04-2024}
			\aligned
			\vert \del^{I_1}L^{J_1}u\vert \vert\del^{I_2}L^{J_2}\del_tv\vert
			&\lesssim (C_1\varepsilon)(s/t)^{1/2}s^{-3/2}(C_1\varepsilon)(s/t)^{3/2}s^{-3/2}s^\delta\\
			&\lesssim(C_1\varepsilon)^2(s/t)^2s^{-3+\delta}.
			\endaligned
		\end{equation}
		For the second term, we remark that in this case $\vert J_2\vert = \vert J\vert - \vert J_1\vert$. Thus
		\begin{equation}\label{eq7-04-04-2024}
			\sum_{{\vert J_1\vert + \vert J_2\vert = \vert J\vert},\atop \vert J_1\vert\geq1}\vert L^{J_1}u\vert\vert\del^{I}L^{J_2}\del_tv\vert\lesssim (s/t)^{2-2\delta}s^{-5/2}\sum_{k_1 = 1}^{k}\Bbf_{k_1}(s)\Abf_{k-k_1}(s).
		\end{equation}
		For the last term, we notice that
		\begin{equation}\label{eq8-04-04-2024}
			\vert [\del^IL^J,\del_t]v\vert\lesssim\vert  v\vert_{N - 4,k-1},
		\end{equation}
		Then we have
		\begin{equation}\label{eq9-04-04-2024}
			\vert u\vert\Big\vert[\del^IL^J,\del_t]v\Big\vert\lesssim (s/t)^{2-2\delta}s^{-5/2}\Bbf_0(s)\Abf_{k-1}(s).
		\end{equation}
		\end{proof}
		
		For the second term of source of \eqref{eq5-24-07-2024}, we establish the following result.
		\begin{lemma}\label{lemma1-24-07-2024}
			We have the following pointwise estimates:
			\begin{equation}\label{eq6-24-07-2024}
				\vert(u_L + u_g)\del_t\del^IL^Jv\vert \lesssim (C_1\varepsilon)^2(s/t)^2s^{-3+\delta}, \quad \vert I\vert + \vert J\vert \leqslant N - 4.
			\end{equation}
		\end{lemma}
		\begin{proof}
		This is based on \eqref{eq9-24-07-2024} together with \eqref{eq5-04-04-2024}:  
			\begin{equation}
				\vert(u_L + u_g)\del_t\del^IL^Jv\vert \lesssim C_1\varepsilon (s/t)^{1/2}s^{-3/2}C_1\varepsilon t^{-3/2}s^\delta = (C_1\varepsilon)^2(s/t)^2s^{-3+\delta}.
			\end{equation}
		\end{proof}
		
		Now we are ready to establish the estimate on $\Abf_k$. 
We have following estimate
	\begin{proposition}\label{Estimate of Abf}
		For $\Abf_k$ and $\Bbf_k$ defined in \eqref{eq1-15-04-2024}, we have following estimates:
		\begin{equation}\label{eq11-20-08-2024}
			\Abf_k(s)\lesssim C_1\varepsilon + \int_{\lambda_0}^{s}\lambda^{-1}\Bbf_0(\lambda)\Abf_{k-1}(\lambda)d\lambda + \sum_{k_1 = 1}^{k}\int_{\lambda_0}^{s}\lambda^{-1}\Bbf_{k_1}(\lambda)\Abf_{k-k_1}(\lambda)d\lambda.
		\end{equation}
		Notice that the last two terms do not exist when $k=0$.
	\end{proposition}
		\begin{proof}
			Remark that by \eqref{eq5-04-04-2024},
			\begin{equation}\label{check01}
				\aligned
				s^{1/2}\vert\del^IL^Jv\vert_1
				\lesssim  s^{1/2}\vert v\vert_{N-3}
				\lesssim  C_1\varepsilon(s/t)^{3/2}s^{-1+\delta}
				\lesssim  C_1\vep (s/t)^{5/2-\delta}.
				\endaligned
			\end{equation}
			Recalling Lemma \ref{Estimate of R1 and R2}, when near the light cone ($r/t\geqslant \frac{s_0^1-1}{s_0^2+1} = \frac{3}{5}$), where $\lambda_0\sim t/s$, we have
			\begin{equation}\label{check02}
				\int_{\lambda_0}^{s}\Big\vert R[u_b,\del^IL^Jv]\Big\vert d\lambda\lesssim C_1\varepsilon(s/t)\int^{s}_{\lambda_0}\lambda^{-2+2\delta}d\lambda\lesssim C_1\varepsilon(s/t)^{2-2\delta}.				
			\end{equation}
			Recall that in the present case, $f: = [\del^IL^J,u\del_t] v + (u_L+ u_g)\del_t \del^IL^J v$. 
To apply proposition~\ref{prop1-04-04-2024}, we combine the results of \eqref{check01}, \eqref{check02}, Lemma~\ref{lemma1-04-04-2024} and Lemma~\ref{lemma1-24-07-2024}
\begin{equation}\label{check04}
\aligned
&s^{3/2}(\vert \del^IL^J v\vert +(s/t)\vert\del\del^IL^J v\vert) 
\\
\lesssim & C_1\varepsilon(s/t)^{5/2-\delta}
+ C_1\varepsilon(s/t)^{2-2\delta}
+ (C_1\vep)^2(s/t)^{5/2-\delta}
\\
& + (s/t)^{2-2\delta}\int_{\lambda_0}^s \lambda^{-1} \Big(\Bbf_0(\lambda)\Abf_{k-1}(\lambda) 
+ \sum^{k}_{k_1 = 1}\Bbf_{k_1}\Abf_{k-k_1}(\lambda)\Big)d\lambda
\\
& + (C_1\vep)^2(s/t)^{5/2-\delta} 
	\endaligned
\end{equation}
Multiplying the above estimate by $(s/t)^{2\delta-2}$, we have 
\begin{equation}\label{eq3-07-03-2025}
	\aligned
	&(s/t)^{2\delta -2}s^{3/2}(\vert \del^IL^J v\vert +(s/t)\vert\del\del^IL^J v\vert) 
	\\
	&\lesssim C_1\vep 
	+ \int_{\lambda_0}^s \lambda^{-1} \Big(\Bbf_0(\lambda)\Abf_{k-1}(\lambda) 
	+ \sum^{k}_{k_1 = 1}\Bbf_{k_1}\Abf_{k-k_1}(\lambda)\Big)d\lambda.
	\endaligned
\end{equation}			
For the case of far light-cone$(0 \leqslant r/t\leqslant \frac{s_0^2-1}{s_0^2+1} = \frac{3}{5}$, $\frac{4}{5}\leqslant s/t\leqslant 1)$, where $\lambda_0 = s_0 = 2$, we notice 
\begin{equation}\label{eq1-07-03-2025}
	\sup_{\mathcal{H}_{s_0}}{(|\del^IL^Jv(t,x)| + |\del \del^IL^Jv(t,x)|)} \lesssim C_1\vep.
\end{equation}
Following the same method, we have 
\begin{equation}\label{eq4-07-03-2025}
	\aligned
	& s^{3/2}(\vert \del^IL^J v\vert +(s/t)\vert\del\del^IL^J v\vert) 
	\\
	\lesssim & C_1\vep
	+ C_1\vep(s/t)^{5/2-\delta}
	+ C\vep (s/t)^{3-\delta}
	+ C_1\vep(s/t)^{3/2} 
	+ (C_1\vep)^2(s/t)
	\\
	& + (C_1\vep)^2(s/t)^{2} 
	+ (C_1\vep)^2(s/t)^{2}
	\\
	& + (s/t)^{2-2\delta}\int_{\lambda_0}^s \lambda^{-1} \Big(\Bbf_0(\lambda)\Abf_{k-1}(\lambda) 
	+ \sum^{k}_{k_1 = 1}\Bbf_{k_1}\Abf_{k-k_1}(\lambda)\Big)d\lambda.
	\endaligned
\end{equation}
That is 
\begin{equation}\label{eq5-07-03-2025}
	\aligned
	&(s/t)^{2\delta -2}s^{3/2}(\vert \del^IL^J v\vert +(s/t)\vert\del\del^IL^J v\vert) 
	\\
	\lesssim & C_1\vep 
	+ \int_{\lambda_0}^s \lambda^{-1} \Big(\Bbf_0(\lambda)\Abf_{k-1}(\lambda) 
	+ \sum^{k}_{k_1 = 1}\Bbf_{k_1}\Abf_{k-k_1}(\lambda)\Big)d\lambda.
	\endaligned
\end{equation}
%
Recall \eqref{eq1-15-04-2024}, the definition of $\Abf_k$, we obtain \eqref{eq11-20-08-2024}.
		\end{proof}

\subsection{The sharp decay estimates for the wave equation}
\begin{proposition}\label{Estimate of Bbf}
	For $\Bbf_k$ defined in \eqref{eq1-15-04-2024}, we have following estimates:
	\begin{equation}\label{eq15-20-08-2024}
		\Bbf_k(s)\lesssim C_1\varepsilon + \sum_{k_1 = 0}^{k}\Abf_{k_1}(s)\Abf_{k-k_1}(s).
	\end{equation}
\end{proposition}
		\begin{proof}
		
Recall our wave equation:
\begin{equation}\label{eq3-21-08-2024}
	\Box (\del^IL^Ju) = \del^IL^J(A^{\alpha\beta}\del_\alpha v\del_\beta v + Bv^2)
\end{equation}
for $\vert I\vert +\vert J\vert \leqslant N - 4$. Here we introduce a decomposition: $\del_IL^Ju = \phi_I + \phi_S$, which satisfy
\begin{equation}\label{eq6-07-03-2025}
	\aligned
	&\begin{cases}
		&\Box \phi_I = 0,
		\\
		&\phi_I\vert_{\mathcal{H}_2} = \del^IL^Ju\vert_{\mathcal{H}_2},
		\\
		&\del_t \phi_I\vert_{\mathcal{H}_2} = \del_t\del^IL^Ju\vert_{\mathcal{H}_2}.
	\end{cases}
	\\
	&\begin{cases}
		&\Box \phi_S = \del^IL^J\big(A^{\alpha\beta}\del_\alpha v\del_\beta v + Bv^2\big),
		\\
		&\phi_S\vert_{\mathcal{H}_2} = 0,
		\\
		&\del_t \phi_S\vert_{\mathcal{H}_2} = 0.
	\end{cases}
	\endaligned
\end{equation}

We remark that $\del^IL^Ju\vert_{\mathcal{H}_2}$ and $\del_t \del^IL^Ju\vert_{\mathcal{H}_2}$ are compactly supported in $\mathcal{H}_2^*$. 
For $\phi_I$ we know that it is the solution to the above free linear wave equation with compactly supported initial data. Thus
\begin{equation}\label{eq2-21-08-2024}
	\vert \phi_I\vert \lesssim C_1\varepsilon t^{-1}.
\end{equation}
		
For $\phi_S$, we want apply Proposition 3.1. For this purpose we need to bound the right-hand-side of \eqref{eq3-21-08-2024}. Thanks to  the definition of $\Abf_k$, for $k\leqslant N-4$
\begin{equation}\label{check06}
	\vert \del v\vert_{N-4,k}\lesssim (s/t)^{1-2\delta}s^{-3/2}\Abf_{k}(s).
\end{equation}
Furthermore, we have 
\\
\begin{equation}\label{check07}
	\vert \del v\del v\vert_{k} \lesssim t^{-5/2+2\delta}(t-r)^{-1/2-2\delta}\sum_{k_1 = 0}^{k}\Abf_{k_1}(s)\Abf_{k-k_1}(s),
\end{equation}
		and
		\begin{equation}\label{check08}
			\vert v ^2\vert_k \lesssim t^{-7/2+2\delta}(t-r)^{1/2-2\delta}\sum_{k_1 = 0}^k\Abf_{k_1}(s)\Abf_{k-k_1}(s).
		\end{equation}
		That is
		\begin{equation}\label{eq4-21-08-2024}
			\vert \del v\del v\vert_{k} + \vert v ^2\vert_k \lesssim t^{-5/2+2\delta}(t-r)^{-1/2-2\delta}
			\sum_{k_1 = 0}^k\Abf_{k_1}(s)\Abf_{k-k_1}(s),\quad in~\mathcal{K}_{[2,s]}.
		\end{equation}
		Then by Proposition \ref{prop-01-10-2023} with $\mu = \nu = 1/2-2\delta$, we obtain
		\begin{equation}\label{eq5-21-08-2024}
			\vert \phi_S\vert\lesssim
			t^{-1}\sum_{k_1 = 0}^k\Abf_{k_1}(s)\Abf_{k-k_1}(s), \quad in~\mathcal{K}_{[2,s]}.
		\end{equation}
		Thus we derive
		\begin{equation}\label{check09}
			\Bbf_k(s)\lesssim C_1\varepsilon + \sum_{k_1 = 0}^{k}\Abf_{k_1}(s)\Abf_{k-k_1}(s).
		\end{equation}
		\end{proof}		
		Furthermore, when $k = 0$, we recall \eqref{eq11-20-08-2024} and \eqref{eq15-20-08-2024}, and obtain
		\begin{equation}\label{eq6-21-08-2024}
			\Abf_0(s) + \Bbf_0(s)\lesssim C_1\varepsilon.
		\end{equation}
		\subsection{Conclusion by induction}
			Substituting \eqref{eq6-21-08-2024}, \eqref{eq11-20-08-2024} and \eqref{eq15-20-08-2024} we derive the following system of integral inequalities:
			\begin{equation}\label{eq7-21-08-2024}
				\aligned
				\Abf_k(s)\leqslant & CC_1\varepsilon +
				CC_1\varepsilon\int_{2}^{s}\lambda^{-1}\Bbf_k(\lambda)d\lambda+ CC_1\varepsilon\int_{2}^{s}\lambda^{-1}\Abf_{k-1}(\lambda)d\lambda 
				\\
				& + C\sum_{k_1 = 1}^{k-1}\int_{2}^{s}\lambda^{-1}\Bbf_{k_1}(\lambda)\Abf_{k-k_1}(\lambda)d\lambda,
				\\
				\Bbf_k(s)\leqslant & CC_1\varepsilon + CC_1\varepsilon\Abf_{k}(s) + C\sum_{k_1 = 1}^{k-1}\Abf_{k_1}(s)\Abf_{k-k_1}(s),
				\endaligned
			\end{equation}
			where $C$ is a constant determined by $N,\delta$. 
			
			By an induction argument, we obtain
			\begin{equation}\label{eq8-21-08-2024}
				\Abf_{k}(s) + \Bbf_k(s)\leqslant CC_1\varepsilon s^{CC_1\varepsilon}, \quad k = 1,2,\cdots,N-4.
			\end{equation}
			
			Here the key structure which permits us to make this induction is as following: the coefficients of the top order terms on the right-hand side (that is, the terms of the same order as those we need to estimate) all exhibit a decay of order $\lambda^{-1}$ (which is sufficient for our decay estimates), while terms of lower order may suffer from a loss of decay. This structure will appear again later and will play a crucial role in the proof of the refined energy estimates. In Appendix~\ref{Induction}, we give a detailed explanation.
			
			Thus we conclude the \eqref{eq11-21-08-2024} and \eqref{eq12-21-08-2024}.

		\section{Refined Energy estimate}
		
		In this section we apply the energy estimates Proposition~\ref{Energy}. For the wave component, we need to establish the following result:
		\begin{equation}\label{eq13-21-08-2024}
			\aligned
			&\Vert \vert \del_\alpha v\del_\beta v\vert_{N,k}\Vert_{L_f^2(\mathcal{H}_s)} + \Vert \vert v^2\vert_{N,k}\Vert_{L_f^2(\mathcal{H}_s)}
			\\
			& \lesssim C_1\varepsilon s^{-3/2}\Ebf_{c}^{N,k}(s,v)^{1/2} + C_1\varepsilon s^{-3/2+CC_1\varepsilon}\Ebf_{c}^{N,k-1}(s,v)^{1/2}.
			\endaligned
		\end{equation}
		This can be checked directly by \eqref{eq11-21-08-2024} and \eqref{eq12-21-08-2024}. We only need to remark that
		\begin{equation}\label{eq14-21-08-2024}
			\aligned
			\vert \del_\alpha v\del_\beta v\vert_{N,k}
			\lesssim & \sum_{p_2+p_2=p\atop k_1+k_2=k}\vert\del v\vert_{p_1,k_1}\vert\del v\vert_{p_2,k_2}
			\\
			\lesssim & \sum_{p_1\leqslant N-4,p_2\geqslant 4\atop k_1+k_2 = k}\vert\del v\vert_{p_1,k_1}\vert \del v\vert_{p_2,k_2}
			\\
			= & \sum_{p_1\leqslant N-4\atop p_2\geqslant 4}\sum_{0\leqslant k_1\leqslant k}\vert \del v\vert_{p_1,k_1}\vert \del v\vert_{p_2,k-k_1}
			\\
			\lesssim  & C_1\vep s^{3/2}(s/t)^{2-2\delta}\vert\del v\vert_{N,k}
			\\
			& + \sum_{1\leqslant k_1\leqslant k}C_1\vep s^{-3/2+CC_1\vep}(s/t)^{2-2\delta}\vert \del v\vert_{N,k-1},
			\endaligned
		\end{equation}
		and then apply \eqref{eq2-21-04-2023}(under the condition $N\geqslant 7$). The estimate of $v^2$ is even easier and we omit the detail. Substitute \eqref{eq13-21-08-2024} into \eqref{Klein-Gordon equation}, and obtain
		\begin{equation}\label{eq8-22-08-2024}
			\aligned
			\Ebf_{\con}^{k}(s,u)^{1/2}
			\leqslant & C_0\vep 
			+ CC_1\vep\int_{s_0}^s\tau\cdot\tau^{-3/2}\Vert(s/t)\vert\del v\vert_{N,k}\Vert_{L^2_f(\Hcal_2)}d\tau
			\\
			& + CC_1\vep\int_{s_0}^s\tau\cdot\tau^{-3/2 + CC_1\vep}\Vert(s/t)\vert\del v\vert_{N,k-1}\Vert_{L_f^2(\Hcal_2)}d\tau
			\\
			= & C_0\varepsilon + CC_1\varepsilon\int_{2}^{s}\tau^{-1/2}\Ebf_c^{N,k}(\tau,v)^{1/2}d\tau
			\\
			& + CC_1\varepsilon\int_{2}^{s}\tau^{-1/2+CC_1\varepsilon}\Ebf_c^{N,k-1}(\tau,v)^{1/2}d\tau.
			\endaligned
		\end{equation}
		When $k=0$, the last term vanishes.
		
		For the Klein-Gordon equation, with respect to $\vert I\vert +\vert J\vert \leqslant N$ and $\vert J\vert = k$, we define
		\begin{equation}\label{eq10-29-08-2024}
			\Box \del^IL^Jv + c^2\del^IL^Jv = \del^IL^J(u\del_tv),
		\end{equation}
		here we notice
		\begin{equation}\label{eq1-01-10-2024}
			\aligned
			 \vert \del^IL^J(u\del_tv) \vert
			 \leqslant & \sum_{{\vert I_1\vert + \vert I_2\vert \leqslant \vert I\vert} \atop\vert J_1\vert +\vert J_2\vert \leqslant \vert J\vert}\vert\del^{I_1}L^{J_1}u\vert \vert \del^{I_2}L^{J_2}\del_tv\vert
			 \\
			 \leqslant & \sum_{\vert I_1\vert +\vert J_1\vert \leqslant N-4\atop \vert I_2\vert + \vert J_2\vert \geqslant 4}\vert \del^{I_1}L^{J_1}u\vert \vert \del^{I_2}L^{J_2}\del_tv\vert
			 \\
			 & + \sum_{\vert I_1\vert +\vert J_1\vert\geqslant 4\atop \vert I_2\vert +\vert J_2\vert \leqslant N-4}\vert \del^{I_1}L^{J_1}u\vert \vert\del^{I_2}L^{J_2}\del_tv\vert,
			 \\
			 : = & T_1 + T_2.
			 \endaligned
		\end{equation}
		For $T_1$ we notice that
		\begin{equation}\label{eq1-07-04-2025}
			T_1 \lesssim \sum_{\substack{\vert I_1\vert +\vert J_1\vert\leqslant N-4\\ \vert I_1\vert \geqslant 1\\ \vert I_2\vert +\vert J_2\vert \geqslant 4}}\vert \del^{I_1}L^{J_1}u\vert \vert \del^{I_2}L^{J_2}\del_tv\vert 
			+ \sum_{\substack{\vert I_1\vert +\vert J_1\vert\leqslant N-4\\ \vert I_1\vert = 0\\ \vert I_2\vert +\vert J_2\vert \geqslant 4}}\vert \del^{I_1}L^{J_1}u\vert \vert \del^{I_2}L^{J_2}\del_tv\vert.
		\end{equation}
		We denote that the first term in the right-hand side as $T_{11}$, so is the second term $T_{12}$. Thanks the \eqref{eq3-04-04-2024}, we have  
		\begin{equation}\label{eq2-07-04-2025}
			\aligned
			T_{11}\lesssim & \sum_{\substack{\vert I_1\vert +\vert J_1\vert\leqslant N-4\\ \vert I_1\vert \geqslant 1\\ \vert I_2\vert +\vert J_2\vert \geqslant 4}}C_1\vep s^{-3/2}(s/t)^{1/2}\vert \del^{I_2}L^{J_2}\del_tv\vert
			\\
			\lesssim & \sum_{\substack{\vert I_1\vert +\vert J_1\vert\leqslant N-4\\ \vert I_1\vert \geqslant 1\\ \vert I_2\vert +\vert J_2\vert \geqslant 4}}C_1\vep s^{-1}(s/t)\vert\del^{I_2}L^{J_2}\del_tv\vert .
			\endaligned
		\end{equation}
		Then we notice that
		\begin{equation}\label{eq3-07-04-2025}
			\Vert T_{11}\Vert_{L^2_f(\Hcal_s)}\lesssim C_1\vep s^{-1}\Ebf_{c}^{N,k}(s,v)^{1/2}.
		\end{equation}
		Similarly, by applying \eqref{eq11-21-08-2024}, we have 
		\begin{equation}\label{eq4-07-04-2025}
			\aligned
			T_{12} = &\sum_{\substack{\vert J_1\vert \leqslant N-4 \\ \vert I\vert +\vert J_2\vert \geqslant 4}} \vert L^{J_1}u\vert \vert \del^{I}L^{J_2}\del_tv\vert 
			\\
			= & \vert u\vert \vert \del^IL^{J}\del_tv\vert 
			+ \sum_{\substack{1\leqslant \vert J_1\vert \leqslant N-4 \\ \vert I\vert +\vert J_2\vert \geqslant 4}}\vert L^{J_1}u\vert \vert\del^IL^{J_2}\del_tv\vert
			\\
			\lesssim & C_1\vep s^{-1}(s/t)\vert \del^IL^{J}\del_tv\vert 
			+ C_1\vep s^{-1+CC_1\vep}(s/t)\sum_{\vert J_2\vert \leqslant \vert J\vert -1}\vert\del^IL^{J_2}\del_tv\vert.
			\endaligned
		\end{equation}
		That is
		\begin{equation}\label{eq5-07-04-2025}
			\aligned
			\Vert T_{12}\Vert_{L^2_f(\Hcal_s)} \lesssim & C_1\vep s^{-1}\Ebf_c^{N,k}(s,v)^{1/2} 
			+ C_1\vep s^{-1 + CC_1\vep}\Ebf_{c}^{N,k-1}(s,v)^{1/2}.
			\endaligned
		\end{equation}
		
		Using the same partitioning method to consider $T_2$, with the result from \eqref{eq5-04-04-2024}, we have
		\begin{equation}\label{eq1-08-04-2025}
			\aligned
			T_{21} = & \sum_{\substack{\vert I_1\vert +\vert J_1\vert\geqslant 4\\ \vert I_1\vert \geqslant 1\\ \vert I_2\vert +\vert J_2\vert \leqslant N-4}}\vert \del^{I_1}L^{J_1}u\vert \vert\del^{I_2}L^{J_2}\del_tv\vert
			\\
			\lesssim & \sum_{\substack{\vert I_1\vert +\vert J_1\vert\geqslant 4\\ \vert I_1\vert \geqslant 1 \\\vert I_2\vert +\vert J_2\vert \leqslant N-4}}C_1\vep s^{-3/2 + \delta}(s/t) \vert\del^{I_1}L^{J_1}u\vert,
			\endaligned
		\end{equation}
		and by applying the result of \eqref{eq12-21-08-2024}, we also have 
		
		\begin{equation}\label{eq3-08-04-2025}
			\aligned
			T_{22} = & \sum_{\substack{\vert J_1\vert \geqslant 4 \\ \vert I\vert  + \vert J_2\vert \leqslant N - 4}}\vert L^{J_1}u\vert\vert \del^{I_2}L^{J_2}\del_tv\vert
			\\
			= & \vert L^{J}u\vert \vert \del^{I}\del_tv\vert 
			+ \sum_{1\leqslant \vert J_2\vert \leqslant N - 4}\vert L^{J_1}u\vert \vert \del^{I}L^{J_2}\del_tv\vert 
			\\
			\lesssim & C_1\varepsilon s^{-3/2}(s/t)\vert L^{J}u\vert 
			+ C_1\varepsilon s^{-3/2+CC_1\varepsilon}(s/t)\sum_{ \vert J_1\vert \leqslant \vert J\vert -1}\vert  L^{J_1}u\vert.
			\endaligned
		\end{equation}
		Then, thanks to \eqref{eq1-14-07-2025}, we have 
		\begin{equation}\label{eq1-10-04-2025}
			\aligned
			\Vert T_{21}\Vert_{L^2_f(\mathcal{H}_s)} + \Vert T_{22}\Vert_{L^2_f(\mathcal{H}_s)} 
			\lesssim & (C_1\vep)^2s^{-3/2 + \delta}
			\\
			& + C_1\varepsilon s^{-3/2}\Ebf_{\con}^{N,k}(s,u)^{1/2}
			+ C_1\varepsilon s^{-3/2+CC_1\varepsilon}\Ebf_{\con}^{N,k-1}(s,u)^{1/2}.
			\endaligned
		\end{equation}
		Finally we have 
		\begin{equation}\label{eq2-10-04-2025}
			\aligned
			\Big\Vert \del^IL^J(u\del_tv) \Big\Vert_{L_f^2(\mathcal{H}_s)}\lesssim
			& (C_1\vep)^2s^{-3/2 + \delta} + C_1\varepsilon s^{-1}\Ebf_{c}^{N,k}(s,v)^{1/2} + C_1\varepsilon s^{-3/2}\Ebf_{\con}^{N,k}(s,u)^{1/2}
			\\
			& + C_1\varepsilon s^{-1+CC_1\varepsilon}\Ebf_{c}^{N,k-1}(s,v)^{1/2}
			\\
			& + C_1\varepsilon s^{-3/2+CC_1\varepsilon}\Ebf_{\con}^{N,k-1}(s,u)^{1/2}.
			\endaligned
		\end{equation}

		Now substitute \eqref{eq2-10-04-2025} into \eqref{eq1-23-08-2023}, we have
		\begin{equation}\label{eq7-22-08-2024}
			\aligned
			\Ebf_c^{N,k}(s,v) \leqslant & C_0\varepsilon + C(C_1\varepsilon)^2 + CC_1\varepsilon\int_{2}^{s}\tau^{-1}\Ebf_{c}^{N,k}(\tau,v)^{1/2}d\tau
			\\
			& + CC_1\varepsilon\int_{2}^{s}\tau^{-3/2}\Ebf_{\con}^{N,k}(\tau,u)^{1/2}d\tau
			\\
			& + CC_1\varepsilon \int_{2}^{s}\tau^{-1 + CC_1\varepsilon}\Ebf_c^{N,k-1}(\tau,v)^{1/2}d\tau
			\\
			& + CC_1\varepsilon\int_{2}^{s}\tau^{-3/2+CC_1\varepsilon}\Ebf_{\con}^{N,k-1}(\tau,u)^{1/2}d\tau.
			\endaligned
		\end{equation}
		When $k=0$, the last term in the above expression vanishes, for the same reason as in \eqref{eq8-22-08-2024}.\\
		The integral inequalities \eqref{eq8-22-08-2024} together with \eqref{eq7-22-08-2024} forms a system. By Grönwall's inequality and proceed by induction on $k$ with $0\leqslant k\leqslant N$, which detailed in the following Appendix~\ref{Induction}, we derive
		\begin{equation}\label{eq9-22-08-2024}
			s^{-1/2}\Ebf_{\con}^{N}(s,u)^{1/2} + \Ebf_c^{N}(s,v)^{1/2}\leqslant 2C_0\varepsilon + C(C_1\varepsilon)^{3/2}s^{C(C_1\varepsilon)^{1/2}},
		\end{equation}
		provided that $C_1>C_0$. If we let
		\begin{equation}\label{eq10-22-08-2024}
			C(C_1\varepsilon)^{1/2}\leqslant \delta,\quad C_1>4C_0,\quad \varepsilon<\frac{(C_1-4C_0)^2}{4C^2C_1^3},
		\end{equation}
		then \eqref{eq9-22-08-2024} leads us to \eqref{Refined energy bound}. That is we proved Proposition~\ref{Bootstrap theorem}.


\begin{appendices}
	\section{Inductive Step}\label{Induction}
	
	To facilitate subsequent calculations, we first provide a preliminary estimate of the formula we intend to compute.
	\begin{equation}
		\aligned
		& s^{-1/2}\Ebf_{\con}^{N,k}(s,u)^{1/2} + \Ebf_c^{N,k}(s,v)^{1/2}\\
		\leqslant & (1+s^{-1/2})C_0\vep + C(C_1\vep)^2 + CC_1\vep\int^s_2\Big(\tau^{-1}\Ebf_{c}^{N,k}(\tau,v)^{1/2} + \tau^{-3/2}\Ebf_{\con}^{N,k}(\tau,u)^{1/2}\Big)d\tau \\
		& + CC_1\vep\int^s_2\Big(\tau^{-1 + CC_1\vep} \Ebf_{c}^{N,k-1}(\tau,v)^{1/2} + \tau^{-3/2 + CC_1\vep}\Ebf_{\con}^{N,k-1}(\tau,u)^{1/2}\Big)d\tau\\
		& + CC_1\vep s^{-1/2}\int_2^s\tau^{-1/2}\Ebf_{c}^{N,k}(\tau,v)^{1/2}d\tau + CC_1\vep s^{-1/2}\int_2^s\tau^{-1/2+CC_1\vep}\Ebf_{c}^{N,k-1}(\tau,v)^{1/2}d\tau\\
		\leqslant & 2C_0\vep + C(C_1\vep)^2 + CC_1\vep \int^s_2\tau^{-1}\Big(\Ebf_{c}^{N,k}(\tau,v)^{1/2} + \tau^{-1/2}\Ebf_{\con}^{N,k}(\tau,u)^{1/2}\Big)d\tau\\
		& + CC_1\vep\int_2^s\tau^{-1 + CC_1\vep}\Big(\Ebf_{c}^{N,k-1}(\tau,v)^{1/2} + \tau^{-1/2}\Ebf_{\con}^{N,k-1}(\tau,u)^{1/2}\Big)d\tau.
		\endaligned
	\end{equation}
	Here we remark that $C_1> C_0$ and in the case of $k = 0$, the last term vanishes. 
	Subsequently, we make induction on $k$ for \eqref{eq9-22-08-2024}. For the special case, $k=0$. We have
	\begin{equation}
		\aligned
		& s^{-1/2}\Ebf_{\con}^{N,0}(s,u)^{1/2} + \Ebf_c^{N,0}(s,v)^{1/2}\\
		\leqslant & 2C_0\vep + C(C_1\vep)^2 + CC_1\vep \int^s_2\tau^{-1}\Big(\Ebf_{c}^{N,0}(\tau,v)^{1/2} + \tau^{-1/2}\Ebf_{\con}^{N,0}(\tau,u)^{1/2}\Big)d\tau,\\
		\endaligned
	\end{equation}
	by applying Grönwall's inequality, we have 
	\begin{equation}
		\aligned
		& s^{-1/2}\Ebf_{\con}^{N,0}(s,u)^{1/2} + \Ebf_c^{N,0}(s,v)^{1/2}\\
		\leqslant & \Big(2C_0\vep + C(C_1\vep)^2\Big)\Big(1 + \int^s_2CC_1\vep\tau^{-1}\exp(\int^s_\tau CC_1\vep \eta^{-1}d\eta)d\tau\Big)\\
		= & \Big(2C_0\vep + C(C_1\vep)^2\Big)\Big(1 + \int^s_2CC_1\vep\tau^{-1}\exp(\ln(s/\tau)^{CC_1\vep})d\tau\Big)\\
		\leqslant & \Big(2C_0\vep + C(C_1\vep)^2\Big)\Big(1 + \int^s_2CC_1\vep\tau^{-1}\exp(\ln(s/t)^{C\sqrt{C_1\vep}})d\tau\Big)\\
		\leqslant & \Big(2C_0\vep + C(C_1\vep)^2\Big)\Big(1 + C(C_1\vep)^{1/2}s^{C\sqrt{C_1\vep}}\Big)\\
		\leqslant & 2C_0\vep + C(C_1\vep)^{3/2}s^{C\sqrt{C_1\vep}}.
		\endaligned
	\end{equation}
	Then we consider the base case $k=1$,
	\begin{equation}
		\aligned
		&s^{-1/2}\Ebf_{\con}^{N,1}(s,u)^{1/2} + \Ebf_c^{N,1}(s,v)^{1/2}\\
		\leqslant & 2C_0\vep + C(C_1\vep)^2 + CC_1\vep \int^s_2\tau^{-1}\Big(\Ebf_{c}^{N,1}(\tau,v)^{1/2} + \tau^{-1/2}\Ebf_{\con}^{N,1}(\tau,u)^{1/2}\Big)d\tau\\
		& + \int^s_2CC_1\vep\tau^{-1+CC_1\vep}\Big(s^{-1/2}\Ebf_{\con}^{N,0}(\tau,u)^{1/2} + \Ebf_c^{N,0}(\tau,v)^{1/2}\Big)d\tau\\
		\leqslant & \Big(2C_0\vep + C(C_1\vep)^2\Big)\Big(1 + \int^s_2CC_1\vep\tau^{-1}\exp(\int^s_\tau CC_1\vep \eta^{-1}d\eta)d\tau\Big)\\
		& +  \int^s_2CC_1\vep\tau^{-1+CC_1\vep}\Big(2C_0\vep + C(C_1\vep)^{3/2}s^{C\sqrt{C_1\vep}}\Big)d\tau\\
		\leqslant & \Big(2C_0\vep + C(C_1\vep)^2\Big)\Big(1 + \int^s_2CC_1\vep\tau^{-1}\exp(\ln(s/t)^{C\sqrt{C_1\vep}})d\tau\Big)\\
		& + \int^s_2CC_1\vep\tau^{-1+C\sqrt{C_1\vep}}\Big(2C_0\vep + C(C_1\vep)^{3/2}s^{C\sqrt{C_1\vep}}\Big)d\tau\\
		\leqslant & \Big(2C_0\vep + C(C_1\vep)^2\Big)\Big(1 + C(C_1\vep)^{1/2}s^{C\sqrt{C_1\vep}}\Big)\\
		& + C(C_1\vep)^{1/2}s^{C\sqrt{C_1\vep}}\Big(CC_1\vep s^{C\sqrt{C_1\vep}}\Big)\\
		\leqslant & \Big(2C_0\vep + C(C_1\vep)^2\Big)\Big(1 + C(C_1\vep)^{1/2}s^{C\sqrt{C_1\vep}}\Big) + C^2(C_1\vep)^{3/2}s^{C^2C_1\vep}\\
		\leqslant & 2C_0\vep + C(C_1\vep)^{3/2}s^{C\sqrt{C_1\vep}}.
		\endaligned
	\end{equation}
	
	Assume that the statement holds for some arbitrary fixed integer $k$, we have following induction hypothesis
	\begin{equation}
		s^{-1/2}\Ebf_{\con}^{N,k}(s,u)^{1/2} + \Ebf_c^{N,k}(s,v)^{1/2}
		\leqslant  2C_0\vep + C(C_1\vep)^{3/2}s^{C\sqrt{C_1\vep}}.
	\end{equation}
	We need to show that the statement then holds for $k+1$,
	\begin{equation}
		\aligned
		&s^{-1/2}\Ebf_{\con}^{N,k+1}(s,u)^{1/2} + \Ebf_c^{N,k+1}(s,v)^{1/2}\\
		\leqslant & 2C_0\vep + C(C_1\vep)^2 + CC_1\vep \int^s_2\tau^{-1}\Big(\Ebf_{c}^{N,k+1}(\tau,v)^{1/2} + \tau^{-1/2}\Ebf_{\con}^{N,k+1}(\tau,u)^{1/2}\Big)d\tau\\
		& + \int^s_2CC_1\vep\tau^{-1+CC_1\vep}\Big(s^{-1/2}\Ebf_{\con}^{N,k}(\tau,u)^{1/2} + \Ebf_c^{N,k}(\tau,v)^{1/2}\Big)d\tau,\\
		\endaligned
	\end{equation}
	by the induction hypothesis and the Grönwall's inequality, we have 
	\begin{equation}
		\aligned
		&s^{-1/2}\Ebf_{\con}^{N,k+1}(s,u)^{1/2} + \Ebf_c^{N,k+1}(s,v)^{1/2}\\
		\leqslant & \Big(2C_0\vep + C(C_1\vep)^2\Big)\Big(1 + \int^s_2CC_1\vep\tau^{-1}\exp(\int^s_\tau CC_1\vep \eta^{-1}d\eta)d\tau\Big)\\
		& +  \int^s_2CC_1\vep\tau^{-1+CC_1\vep}\Big(2C_0\vep + C(C_1\vep)^{3/2}s^{C\sqrt{C_1\vep}}\Big)d\tau\\
		\leqslant & \Big(2C_0\vep + C(C_1\vep)^2\Big)\Big(1 + \int^s_2CC_1\vep\tau^{-1}\exp(\ln(s/t)^{C\sqrt{C_1\vep}})d\tau\Big)\\
		& + \int^s_2CC_1\vep\tau^{-1+C\sqrt{C_1\vep}}\Big(2C_0\vep + C(C_1\vep)^{3/2}s^{C\sqrt{C_1\vep}}\Big)d\tau\\
		\leqslant & \Big(2C_0\vep + C(C_1\vep)^2\Big)\Big(1 + C(C_1\vep)^{1/2}s^{C\sqrt{C_1\vep}}\Big) + C^2(C_1\vep)^{3/2}s^{C^2C_1\vep}\\
		\leqslant & 2C_0\vep + C(C_1\vep)^{3/2}s^{C\sqrt{C_1\vep}}.
		\endaligned
	\end{equation}
	Hence, by induction, the result follows for all $0\leqslant k\leqslant N$.
\end{appendices}

\addcontentsline{toc}{section}{References}

\pagestyle{plain}

\end{document}